\documentclass[12pt]{article}
\usepackage{amsfonts}
\usepackage{amssymb}
\usepackage{amsmath, amsthm}

\makeatletter
\@addtoreset{equation}{section}
\makeatother

\newtheorem{theorem}{Theorem}[section]
\newtheorem{lemma}[theorem]{Lemma}
\newtheorem{corollary}[theorem]{Corollary}
\newtheorem{proposition}[theorem]{Proposition}

\theoremstyle{definition}
\newtheorem{definition}[theorem]{Definition}
\newtheorem{example}[theorem]{Example}

\newtheorem{remark}[theorem]{Remark}

\newtheorem*{pro_f}{\it Proof}

\begin{document}

\noindent{\LARGE \textbf{Some Combinatorial Aspects of
 \medskip }}

\noindent {\LARGE \textbf{Constructing Bipartite-Graph
Codes}\bigskip}

\noindent \textbf{Alexander A. Davydov,$^{\textbf 1}$ Massimo
Giulietti,$^{\textbf 2}$
 Stefano Marcugini,$^{\textbf
2}$\newline Fernanda Pambianco$^{\textbf 2}$}\newline
$^{1}$\textsl{Institute for Information Transmission Problems,
Russian Academy of Sciences, Bol'shoi Karetnyi per. 19, GSP-4,
Moscow, 127994, Russian Federation, E-mail: adav@iitp.ru
\medskip }

\noindent$^{2}$\textsl{Dipartimento di Matematica e
Informatica, Universit\`{a} degli Studi di Perugia, Via
Vanvitelli 1, Perugia, 06123, Italy
 \newline E-mail:
giuliet@dipmat.unipg.it; gino@dipmat.unipg.it;
fernanda@dipmat.unipg.it }\medskip

\noindent\textbf {Abstract: We propose geometrical methods for
constructing square 01-matrices with the same number $n$ of
units in every row and column, and such that any two rows of
the matrix contain at most one unit in common. These matrices
are equivalent to $n$-regular bipartite graphs without
4-cycles, and therefore can be used for the construction of
efficient bipartite-graph codes such that both the classes of
its vertices are associated with local constraints. We
significantly extend the region of parameters $m,n$ for which
there exist an $n$-regular bipartite graph with $2m$ vertices
and without $4$-cycles. In that way we essentially increase the
region of lengths and rates of the corresponding
bipartite-graph codes. Many new matrices are either circulant
or consist of circulant submatrices: this provides code
parity-check matrices consisting of circulant submatrices, and
hence quasi-cyclic bipartite-graph codes with simple
implementation.}\bigskip


\section{Introduction}

Bipartite-graph codes are studied in the context of low-density parity-check
(LDPC) codes, i.e., error correcting codes with a strongly sparse parity
check matrix. These codes were first presented by Gallager \cite{Galager} in
1962, see also \cite{Tanner81}-\cite{DGMP-ACCT2008} and the references
therein.

The general idea of connecting linear LDPC codes to bipartite
graphs first appeared in Tanner's seminal paper
\cite{Tanner81}. In the Tanner graph of an LDPC code
$\mathcal{C}$, the vertices of one class (variable vertices)
correspond to code symbols and those of the other class are
associated with subcodes (local constraints on variables). The
case where the subcodes are linear codes with a single
parity-check has been intensively investigated, see e.g.
\cite{Tanner81},\cite{KouLinFos-RediscNew},\cite{JonWel-Resolv2-des}
-\cite{CodFinGeom2005},\cite{GabidISIT}-\cite{QCFiniteField}.
In this case, each parity-check of the code $\mathcal{C}$ is
represented by a subcode vertex.

In a number of papers, see e.g.
\cite{SipSp}-\cite{BPZem},\cite{MilFos},\cite
{MilFos2008},\cite{ZLJR-1}, codes $\mathcal{C}$ in which the
local constraints on variables have few parity-checks are
considered. The subcodes are called component or constituent
codes, and $\mathcal{C}$ is said to be a generalized LDPC
(GLDPC) code. It is natural to refer to the corresponding
bipartite graph as to the \emph{generalized Tanner graph} of
$\mathcal{C}$.

Sipser and Spielman \cite{SipSp} significantly developed the approach in
\cite{Tanner81} by using graph expansion parameters and spectral properties
of the Tanner graph for decoding analysis; they also suggested the term
\emph{expander codes} for code families whose analysis relies on graph
expansion.

Variants of Tanner's construction were proposed in several
papers. In the codes $\mathcal{C}$ investigated in
\cite{BargZemor2002},\cite{BargZemor},\cite{HohJust},\cite{AfDaZ}-\cite{BargMazu},
the code symbols correspond to the edges of a bipartite graph
$G$; each vertex $v$ of $G$ is associated with a local
constraint $\mathcal{C}_{v}$ of length equal to the degree of
the vertex. The subcode $\mathcal{C}_{v}$ coincides with the
projection of $ \mathcal{C}$ on the positions corresponding to
the edges containing $v$. Usually, the graph $G$ is assumed to
be regular (i.e. each vertex has the same degree) or biregular
(i.e. each vertex from the same class has the same degree).
This allows the use of the same local constraint for each
vertex $v$ from the same class, and therefore facilitates
decoding analysis. According to \cite{BargZemor} and
\cite{BargMazu}, such a code is called a \emph{ bipartite-graph
code }(BG code for short). It should be noted that the graph
$G$ is distinct from the generalized Tanner graph of
$\mathcal{C}$. Throughout the paper we refer to $G$ as the
\emph{supporting graph} of $ \mathcal{C}$.

In order to improve the performance of a bipartite-graph code,
it is desirable to increase the \emph{girth} of the graph, that
is, the minimum length of its cycles
\cite{KouLinFos-RediscNew},\cite{JonWel-Resolv2-des}-\cite
{CodFinGeom2005},\cite{GabidISIT}-\cite{AfDaZ}.
In this paper, the combinatorial aspects of this issue are
dealt with. In particular, we investigate the spectrum $\Sigma
$ of parameters $m_{1},m_{2},n_{1},n_{2}$ for which there
exists a biregular bipartite graph $(V_{1}\cup V_{2},E)$ free
from $4$-cycles and such that $|V_{i}|=m_{i}$, $\deg (v)=n_{i}$
for all $v$ in $V_{i}$. A number of new geometrical methods for
constructing these graphs are proposed, so that the spectrum of
parameters of bipartite-graph codes is significantly extended,
especially for the case $m_{1}=m_{2}$, $ n_{1}=n_{2}$. It
should be noted that if the supporting graph of a BG code has
girth at least\emph{\ }six, then the girth of the generalized
Tanner graph of this code (considered as a GLDPC code) is at
least ten.

A key tool in our investigation is the $01$-matrix $M(G)$
corresponding to a biregular bipartite graph $G$ with
parameters $m_{1},m_{2},n_{1},n_{2}$ and without $4$-cycles:
each row corresponds to a vertex in $V_{1}$, each column to a
vertex in $V_{2}$, and the entry in position $(i,j)$ is $1$ if
and only if there is an edge joining the vertices corresponding
to the $i$-th row and the $j$-th column. The matrix $M(G)$ is
an $m_{1}\times m_{2}$ matrix, with $ n_{1}$ units in every row
and $n_{2}$ units in every column; also, since $G$ has no
$4$-cycles, the $2\times 2$ matrix $J_{4}$ consisting of all
units is not a submatrix of $M(G)$. We denote the class of such
matrices as $ M(m_{1},m_{2},n_{1},n_{2})$. For the sake of
simplicity, we write $M(m,n)$ for $M(m,m,n,n)$. Clearly, from a
matrix $M$ in $M(m_{1},m_{2},n_{1},n_{2})$ one can construct a
bipartite biregular graph $G(M)$ free from $4$-cycles with
parameters $m_{1},m_{2},n_{1},n_{2}.$ Most of the paper is
devoted to the construction of matrices of type
$M(m_{1},m_{2},n_{1},n_{2})$ - and hence of bipartite biregular
graphs with parameters $m_{1},m_{2},n_{1},n_{2}$ - mainly based
on incidence structures in finite projective spaces $PG(v,q)$
over Galois fields $F_{q}$ (see \cite{Hirs},\cite{Handbook} for
basic facts on Galois Geometries).

We remark that $J_{4}$-free matrices are studied in connection
with LDPC codes not just as matrices of supporting graphs, see
e.g. \cite
{KouLinFos-RediscNew},\cite{JonWel-Resolv2-des}-\cite{CodFinGeom2005},\cite
{GabidISIT}-\cite{AfDaZ}, and the references therein. Mainly,
non-square matrices are investigated; exceptions can be found
in \cite{KPPPF-IEEE},\cite
{GabidISIT}-\cite{GabZvenig},\cite{AfDaZ}. Matrices of type
$M(m,n)$ are also relevant in design theory. The fact that the
incidence matrix of a $2$-$ (v,k,1)$ design which is either
symmetric or resolvable non-symmetric (see \cite{Handbook})
gives rise to matrices of type $M(v,k)$ is well known, and has
already been used in some works on LDPC codes, see \cite
{JonWel-Resolv2-des, OvalDes-WelJon, JW-Unital}. We point out
that weaker incidence structures such as \emph{symmetric
configurations} (see \cite[ Section
IV.6]{Handbook},\cite{Gropp-nk}-\cite{FunkLabNap}) give rise to
matrices of type $M(m,n)$ as well, and we will significantly
rely on this remark throughout the paper. Also, methods from
Graph Theory have turned out to be very useful for constructing
of matrices of type $M(m,n)$, see e.g.
\cite{KPPPF-IEEE},\cite{AFLN-graphs},\cite{GH}, and the
references therein.

Despite matrices $M(m_1,m_2,n_1,n_2)$ having been thoroughly investigated,
the spectrum of parameters in $\Sigma$ seems to be not wide enough if
compared to the permanently growing needs of practice, where exact values of
$m_i$ and $n_i$ are often necessary. Also, it should be taken into
consideration that distinct constructions provide matrices with distinct
properties, and clearly some choice can be useful.

The paper is organized as follows. First, we provide a construction
(Construction A) based on incidence structures whose point-set consists of a
single orbit of points under the action of a collineation group of a finite
affine or projective plane, see Section \ref{Sec3}.

In Section \ref{Sec4} we provide Construction B. The starting
point is any cyclic symmetric configuration $\mathcal{I}$. Then
we consider incidence substructures in which the point-set is
the union of orbits under the action of subgroups of the
automorphism group of $\mathcal{I}$. Construction~B turns out
to be very productive, as it gives rise to numerous matrices $
M(m,n)$ and $M(m_{1},m_{2},n_{1},n_{2})$ with distinct
parameters. Moreover, the new matrices $M(m,n)$ are either
\emph{circulant }or consist of \emph{ circulant submatrices}.
Therefore, not only the spectrum $\Sigma $ is significantly
extended, but also parity-check matrices with interesting and
useful structure are provided, see Remarks
\ref{Rem4_circ-weight-h} and \ref {Rem4_skeleton_circulant}. In
particular, we obtain parity-check matrices consisting of
circulant submatrices. The corresponding bipartite-graph codes
are \emph{quasi-cyclic }(QC). Interestingly, QC codes can be
encoded with the help of shift-registers with complexity
linearly proportional to their code length, see
\cite{LinCost},\cite{QCEncoder}.

In Sections \ref{Sec5} and \ref{Sec5.5} we apply the general
results obtained in Section \ref{Sec4} to two well-known
classes of cyclic symmetric configurations. This allows us to
obtain further parameters of matrices $ M(m,n)$ and
$M(m_{1},m_{2},n_{1},n_{2}).$ We give also a number of computer
results.

Finally, in Section \ref{Sec6} constructions not arising from collineation
groups are described, and in Section \ref{Sec8} a summary of new symmetric
configurations is given.

Some of the results from this work were presented without
proofs in \cite {DGMP-ACCT2008}.

\section{\label{Sec2}Preliminaries. Problem Statement}

\subsection{\emph{\ }Constructions of Codes}

We consider some constructions of codes on bipartite graphs. Let an $[n,k]$
code be a linear code of length $n$ and dimension $k.$ Note that the methods
developed in the present paper are suitable for both binary and non-binary
codes. \medskip

\emph{Tanner graphs.} In \cite{Tanner81} Tanner proposed the
following way to associate a bipartite graph $T$ to an $[N,K]$
code $\mathcal{C}$ defined by an $R\times N$ parity-check
matrix $H$ with $R\geq N-K.$ One class of vertices
$\{V_{1}^{\prime },\ldots ,V_{N}^{\prime }\}$ of $T,$ called
\emph{ variable vertices}, corresponds to the $N$ positions of
the codewords of$~{ \mathcal{C}}$. Every vertex $V_{i}^{\prime
\prime }$ of the other class $ \{V_{1}^{\prime \prime
},V_{2}^{\prime \prime },\ldots ,V_{R}^{\prime \prime }\}$ is
called a \emph{subcode vertex }and is associated with the
parity-check relation corresponding to the $i$-th row of $H.$
In other words, the $j$-th column ($i$-th row) of $H$ is
identified with a vertex $ V_{j}^{\prime }$ ($V_{i}^{\prime
\prime }$) and a nonzero entry on a position $(i,j)$ of $H$
determines an edge of$~T$ connecting $V_{j}^{\prime } $
and$~V_{i}^{\prime \prime }$. Note that a subcode corresponding
to a vertex $V_{i}^{\prime \prime }$ can be treated as an
$[n_{i},n_{i}-1]$ code $ \mathcal{C}_{i}$ of codimension one
whose parity-check matrix consists of the $i$-th row of $H.$
This construction is considered in numerous works, see e.g.
\cite{KouLinFos-RediscNew},\cite{JonWel-Resolv2-des}-\cite
{CodFinGeom2005},\cite{GabidISIT}-\cite{QCFiniteField}, and the
references therein. The graph $T$ is called a \emph{Tanner
graph} of the code$~\mathcal{ C}$.

The above construction was generalized as follows, see e.g.
\cite{SipSp}-\cite{BPZem},\cite{MilFos},\cite{MilFos2008},\cite{ZLJR-1}.
To a vertex $ V_{i}^{\prime \prime }$ it is associated an
$[n_{i},k_{i}]$ subcode $ \mathcal{C}_{i}$ with $n_{i}-k_{i}>1$
redundancy symbols. In this case subcodes are called component
or constituent codes and the code $\mathcal{C}$ is said to be a
\emph{generalized }LDPC\ code. We refer to the bipartite graph
$T$ as to the \emph{generalized Tanner graph} of $\mathcal{C}$.
The degree of $V_{i}^{\prime \prime }$ is equal to the length
$n_{i}$, and the code $\mathcal{C}_{i}$ is the projection of
$\mathcal{C}$ on the positions associated with the vertices
$V_{j_{1}}^{\prime },\ldots ,V_{j_{n_{i}}}^{\prime }$ adjacent
to~$V_{i}^{\prime \prime }$. Let $R^{\ast }$ denote the number
of constituent codes. Clearly, $R=\sum_{i=1}^{R^{\ast
}}(n_{i}-k_{i})$ holds. The binary $R^{\ast }\times N$ matrix
$S$, in which columns (rows) correspond to vertices
$V_{j}^{\prime }$ ($V_{i}^{\prime \prime }$) and the entry in
position $(i,j)$ is $1$ if and only if there is an edge of $T$
joining $V_{j}^{\prime }$ and $V_{i}^{\prime \prime }$, will be
called a \emph{skeleton }(\emph{framework}) for the $R\times N$
parity-check matrix~$H$. Note that $H$ can be obtained from $S$
by substituting every unit on the $i$-th row of $S$ with a
($(n_{i}-k_{i})$ -positional) column of a parity-check matrix
$H_{i}$ of $\mathcal{C}_{i}.$

Another connection between codes and bipartite graphs is now
considered, see e.g.
\cite{BargZemor2002},\cite{BargZemor},\cite{HohJust},\cite{AfDaZ}-\cite
{BargMazu}. The main difference with respect to Tanner's
construction is that a subcode is associated to \emph{each}
vertex of the graph, whereas code-symbols correspond to its
edges. \medskip

\emph{Basic construction of BG codes} (\cite{BargZemor}). Let
$G$ be an $n$ -regular bipartite graph with two classes of
vertices $\{V_{1},\ldots ,V_{m}\}$ and $\{V_{m+1},\ldots
,V_{2m}\}$ (i.e. any vertex is adjacent to exactly $n$
vertices, but any two vertices from the same class are not
adjacent). Let $\mathcal{C}_{t}$ be an $[n,k_{t}]$
\emph{constituent code}, $ t=1,2,\ldots ,2m.$ A
\emph{bipartite-graph code} $\mathcal{C=C}(G;$ $
\mathcal{C}_{1},\ldots ,\mathcal{C}_{2m})$ is a linear $[N,K]$
code with length equal to the number of edges of $G$, that is
$N=mn$. Coordinates of $ \mathcal{C}$ are in one-to-one
correspondence with the edges of $G$. In addition, the
projection of a codeword of ${\mathcal{C}}$ on the positions
corresponding to the $n$ edges incident to the vertex $V_{t}$
must be a codeword of the constituent code ${\mathcal{C}}_{t}$.
We call $G$ a \emph{\ supporting graph }of the bipartite-graph
code~${\mathcal{C}}$. Straightforward relations between the
parameters of $G$ and $\mathcal{C}$ are the following:
\begin{equation}
N=mn,\text{ }K\leq N-\sum_{t=1}^{2m}(n-k_{t})=\sum_{t=1}^{2m}k_{t}-mn.
\label{form2_NK}
\end{equation}

In \cite{HohJust} the case $\mathcal{C}_{1}=\ldots
=\mathcal{C}_{2m}$ is considered. In general, the constituent
codes can be different. In \cite {BargZemor} the case
$\mathcal{C}_{1}=\ldots =\mathcal{C}_{m},$ $\mathcal{C}
_{m+1}=\ldots =\mathcal{C}_{2m},$ $k_{1}\neq k_{m+1},$ is
investigated. In \cite{AfDaZ} $k_{1}=\ldots =k_{2m},$ but $2m$
distinct generalized Reed-Solomon codes $\mathcal{C}_{t}$ are
used. \medskip

\emph{Generalized basic construction of BG codes} (see e.g.
\cite{BargMazu} ). Let the supporting bipartite graph $G$ be
biregular, and let $ \{V_{1},\ldots ,V_{m_{1}}\}$ and
$\{V_{m_{1}+1},\ldots ,V_{m_{1+}m_{2}}\}$ be its classes of
vertices. Also, let the degree of every vertex from the first
(second, resp.) class be equal to $n_{1}$ ($n_{2}$, resp.).
Then $ m_{1}n_{1}=m_{2}n_{2}.$ The constituent code
$\mathcal{C}_{t}$ is an $ [n_{1},k_{t}]$ code for $t\leq m_{1}$
and an $[n_{2},k_{t}]$ code for $ t>m_{1}.$ The parameters of
the $[N,K]$ BG code $\mathcal{C}$ satisfy
\begin{equation*}
N=m_{1}n_{1}=m_{2}n_{2},
\end{equation*}
\begin{equation}
K\leq
N-\sum_{t=1}^{m}(n_{1}-k_{t})-\sum_{t=m+1}^{2m}(n_{2}-k_{t})=
\sum_{t=1}^{2m}k_{t}-m_{1}n_{1}.  \label{form2_N K variant2}
\end{equation}

\subsection{Parity-check matrices of BG codes}

The parity-check matrix $H$ of a BG code $\mathcal{C}$ can be
obtained in two steps, see \cite{AfDaZ}. First, a skeleton
matrix of $ H$ is constructed. From codes arising from the
basic construction, such skeleton matrix is a binary $2m\times
N$ matrix $S(m,n)$ in which the $t$-th row is associated with
the vertex $V_{t}$. The $j$-th column contains two units in the
positions corresponding to the vertexes incident with the $j$
-th edge. In other words, $S(m,n)$ is the incidence matrix
(\textquotedblleft \emph{edges-vertices}\textquotedblright ) of
the supporting graph $G(M(m,n)).$ Then, $H$ can be obtained
from $S(m,n)$ by substituting every unit on the $t$-th row of
$S(m,n)$ with a ($(n-k_{t})$ -positional) column of a
parity-check matrix $H_{t}$ of $\mathcal{C}_{t}.$ For codes
arising from the generalized basic construction the procedure
is similar.

In Sections \ref{Sec4},\ref{Sec5} we construct \emph{circulant
}matrices $ M(m,n).$ The corresponding skeleton $S(m,n)$ can be
represented as a block matrix consisting of $n$ identity
matrices of order $m$ and $n$ circulant permutation matrices of
size $m\times m$, see Remarks \ref {Rem4_circ-weight-h} and
\ref{Rem4_skeleton_circulant}. Moreover, if the constituent
codes are such that $\mathcal{C}_{1}=\ldots =\mathcal{C}_{m}$
and $\mathcal{C}_{m+1}=\ldots = \mathcal{C}_{2m},$ then the
parity-check matrix $H$ obtained from $S(m,n)$ consists of
circulant submatrices and defines a QC code with relatively
simple implementation.

From the skeleton $S$ of a BG code (considered as a GLDPC code) it is
straightforward to obtain its generalized Tanner graph. In fact, it is the
graph whose adjacency matrix is$~S$.

\subsection{Problem Statement}

Let $M$ be a matrix of type either $M(m,n)$\emph{\ }or $
M(m_{1},m_{2},n_{1},n_{2})$, and let $G=G(M)$. Since $M$ is
$J_{4}$-free, the graph $G$ has girth at least $6.$ This
improves performance of BG code $ \mathcal{C}$ whose supporting
graph is $G$. Also, if $\mathcal{C}$ is viewed as a GLDPC code,
it is easily seen that the girth of its generalized Tanner
graph is at least~$10$. By (\ref{form2_NK}),(\ref{form2_N K
variant2}), parameters $N$ and $K$ of the BG code $\mathcal{C}$
depend on the values of $ m,n,m_{1},m_{2},n_{1},n_{2}.$

Our \emph{goal }is to provide \emph{as many distinct parameters of BG codes
with 4-cycle-free supporting graphs as possible}. This is equivalent to
construct $J_{4}$-free matrices $M(m,n)$ and $M(m_{1},m_{2},n_{1},n_{2})$
with as many distinct parameters as possible.

\subsection{Incidence structures}

Many of the new matrices of type $M(m_1,m_2,n_1,n_2)$
constructed in this paper arise from incidence matrices of
incidence structures. An incidence structure is a pair
$\mathcal{I}=(\mathcal{P},\mathcal{L})$, where $\mathcal{ \ P}$
is a set whose elements are called \emph{points}, and
$\mathcal{L}$ is a collection of subsets of $\mathcal{P}$
called \emph{blocks} (or \emph{lines }). A point $P$ and a
block $\ell$ are said to be \emph{incident} if $P\in \ell$. An
incidence matrix of $\mathcal{I}$ is a $01$-matrix
$M(\mathcal{I})$ where rows corresponds to blocks, columns to
points, and an entry is $1$ if and only if the corresponding
point belongs to the corresponding block.

It is easily seen that the matrix $M(\mathcal{I})$ is a matrix
of type $ M(m_{1},m_{2},n_{1},n_{2})$ if and only if the
incidence structure $\mathcal{ \ I}$ satisfies the following
properties: the number of points in each block is a constant
$n_{1}$, the number of blocks containing a point is a constant
$n_{2}$, and no point pair is contained in more than one block
(note that $ m_{1}$ is the number of blocks, and $m_{2}$ is the
number of points). An incidence structure with this property is
said to be a \emph{configuration } \cite[Section
IV.6]{Handbook}.

If $n_1=n_2=k$ (or, equivalently, the configuration has the
same number $v$ of points and blocks), then $\mathcal{I}$ is
said to be a \emph{symmetric $ (v,k)$-configuration}. In this
case, $M(\mathcal{I})$ is of type $M(m,n)$, where
$m=|\mathcal{P}|=|\mathcal{L}|$ and $n=n_1=n_2$. Projective and
affine spaces over finite fields are well-known examples of
symmetric configurations.

\section{\label{Sec3}A Geometrical Construction of $J_{4}$-free matrices}

Incidence matrices of several geometrical structures are widely
used for obtaining parity-check matrices of LDPC codes, see
\cite{KouLinFos-RediscNew},\cite{JonWel-Resolv2-des}-\cite{JW-Unital},\cite{Fos-circul},\cite
{CodFinGeom2005},\cite{HohJust}-\cite{DGMP-ACCT2008} and the
references therein.

In this section, we consider some configurations arising from
projective and affine spaces, in order to obtain matrices of
types $M(m,n)$ and $ M(m_1,m_2,n_1,n_2)$. \medskip

\textbf{Construction A}. Take any point orbit $\mathcal{P}$
under the action of a collineation group $\Gamma $ in an affine
or projective space of order $ q$. Choose an integer $n\leq
q+1$ such that the set $\mathcal{L}(\mathcal{P} ,n)$ of lines
meeting $\mathcal{P}$ in precisely $n$ points is not empty.
Define the following incidence structure $\mathcal{I}$: the
points are the points of $\mathcal{P}$, the lines are the lines
of $\mathcal{L}(\mathcal{P},n)$,  the incidence is that of the
ambient space. Let $M$ be the incidence matrix of
$\mathcal{I}$.

\begin{theorem}
\label{Th3_constrOneOrb} The incidence structure $\mathcal{I}$
in Construction A is a configuration such that its incidence
matrix is a matrix $M(|\mathcal{P}|,n)$ or $M(|\mathcal{L}(
\mathcal{P},n)|,|\mathcal{P} |,n,r_{n})$, where $r_n$ denotes
the number of lines of $\mathcal{L}( \mathcal{P},n)$ through a
point in $\mathcal{P}$.
\end{theorem}

\begin{pro_f}
The number of points in a block is $n$. The number of lines
through a point is a constant $r_{n}$: this depends on the key
property that the collineation group $\Gamma $ acts
transitively on $\mathcal{P}$. Finally, no two points lie in
more than one block. We obtain a matrix $M(|\mathcal{P} |,n) $\
if $n=r_{n}$ or a matrix
$M(|\mathcal{L}(\mathcal{P},n)|,|\mathcal{P} |,n,r_{n})$
otherwise.$\hfill \qed$
\end{pro_f}

\begin{example}
We consider a regular hyperoval $\mathcal{O}$ in the projective
plane $ PG(2,q)$, $q$ even, see \cite[Section 8.4]{Hirs}. Let
$\mathcal{P} =PG(2,q)\setminus \mathcal{O}$. Then
$|\mathcal{P}|=q^{2}-1.$ The set $ \mathcal{P}$ is an orbit
under the action of the collineation group $\Gamma \cong
PGL(2,q)$ fixing $\mathcal{O}$. Let $n=q+1.$ The set
$\mathcal{L}( \mathcal{P},q+1)$ consists of $\frac{1}{2}q(q-1)$
lines external to $ \mathcal{O}.$ Every point of $\mathcal{P}$
lies on $\frac{1}{2}q$ such lines, as $\mathcal{O}$ has no
tangents. We obtain a matrix of type
\begin{equation*}
M(m_{1},m_{2},n_{1},n_{2}) :m_{1}=\frac{q(q-1)}{2},\text{ }
m_{2}=q^{2}-1,\text{ } n_{1}=q+1,\text{ }n_{2}=\frac{q}{2},\text{ }q\text{ even.}
\end{equation*}
Note that the dual incidence structure is the
$2$-$(\frac{1}{2}q(q-1),$ $ \frac{1}{2}q,1)$ oval design
\cite{OvalDes-WelJon}.
\end{example}

\begin{example}
We consider a conic $\mathcal{K}$ in $PG(2,q)$, $q$ odd, see \cite[Section
8.2]{Hirs}.

\begin{enumerate}
\item Let $\mathcal{P}$ be the set of $\frac{1}{2}q(q-1)$ \emph{internal
points }to $\mathcal{K}$. It is an orbit under the action of the
collineation group $\Gamma _{\mathcal{K}}\cong PGL(2,q)$ fixing the conic.

For $n=\frac{1}{2}(q-1),$ the set
$\mathcal{L}(\mathcal{P},n)$ consists of $
\frac{1}{2}q(q+1)$ bisecants of $\mathcal{K}.$ Every
internal point lies on $ r_{n}=\frac{1}{2}(q+1)$ bisecants.
We obtain a matrix of type
\begin{eqnarray*}
M(m_{1},m_{2},n_{1},n_{2}):m_{1}=\frac{q(q+1)}{2},\,m_{2}=\frac{
q(q-1)}{2},\,n_{1}=\frac{q-1}{2},\,n_{2}=\frac{q+1}{2},\,q\text{ odd.}
\end{eqnarray*}
For $n=\frac{1}{2}(q+1),$ we take as $\mathcal{L}(\mathcal{P},n)$ the set of
$\frac{1}{2}q(q-1)$ lines external to $\mathcal{K}.$ Every internal point
lies on $r_{n}=\frac{1}{2}(q+1)$ external lines. We obtain a matrix of type
\begin{equation*}
M(m,n):m=\frac{q(q-1)}{2},\text{ }n=\frac{q+1}{2},\text{ }q\text{ odd.}
\end{equation*}

\item Another orbit $\mathcal{P}_{2}$ of the group $\Gamma
    _{\mathcal{K}}$ is the set of $\frac{1}{2}q(q+1)$
    \emph{external points} to$~\mathcal{K}$. Let
    $n=\frac{1}{2}(q-1).$ We form the set
    $\mathcal{L}(\mathcal{P}_{2},\frac{ 1}{2}(q-1))$ from
    $\frac{1}{2}q(q+1)$ bisecants. Every external point
    lies on $r_{n}=n=\frac{1}{2}(q-1)$ bisecants. We obtain
    a matrix of type
\begin{equation*}
M(m,n):m=\frac{q(q+1)}{2},\text{ }n=\frac{q-1}{2},\text{ }q\text{ odd.}
\end{equation*}
\end{enumerate}
\end{example}

\begin{example}
\label{Ex3_Hermit} In $PG(2,q)$, with $q$ a square, let
$\mathcal{P}$ be the complement of the Hermitian curve
$\mathcal{H}$ of equation $x_{0}^{\sqrt{q}
+1}+x_{1}^{\sqrt{q}+1}+x_{2}^{\sqrt{q}+1}=0$ \cite[Section
7.3]{Hirs},\cite[ Section 7.11]{Handbook}. It is an orbit under
the action of the projective unitary group $PGU(3,q)$. As
$\mathcal{H}$ contains $q\sqrt{q}+1$ rational points, we have
$|\mathcal{P}|=q^{2}+q-q\sqrt{q}.$ In $PG(2,q)$ there are $
q^{2}+q-q\sqrt{q}$ lines meeting $\mathcal{H}$ in $\sqrt{q}+1$
points and $ \mathcal{P}$ in $q-\sqrt{q}$ points. Every point
of $\mathcal{P}$ lies on $q- \sqrt{q}$ such lines. We obtain a
matrix of type
\begin{equation*}
M(m,n):m=q^{2}+q-q\sqrt{q},\text{ }n=q-\sqrt{q},\text{ }q\text{ square.}
\end{equation*}
The remaining $q\sqrt{q}+1$ lines of $PG(2,q)$ are tangent to
$\mathcal{H}$ and meet $\mathcal{P}$ in $q$ points. Every point
of $\mathcal{P}$ lies on $ \sqrt{q}+1$ such lines. We obtain a
matrix of type
\begin{equation*}
M(m_{1},m_{2},n_{1},n_{2}):m_{1}=q\sqrt{q}+1,\text{ }m_{2}=q^{2}+q-q
\sqrt{q},\,\, n_{1}=q,\text{ }n_{2}=\sqrt{q}+1,\text{ }q\text{ square.}
\end{equation*}
\end{example}

It should be noted that Construction A works for any $2$-$(v,k,1)$ design $D$
and for any group of automorphisms of $D$. The role of $q+1$ is played by
the size of any block in$~D$.

\section{\label{Sec4}A construction from cyclic symmetric configurations}

The aim of this section is to use automorphisms of symmetric
$(v,k)$ -configurations $\mathcal{I}$ in order to obtain
matrices of type $M(m,n)$ other than $M(\mathcal{I})$. In
particular, we are interested in cyclic automorphism groups of
$\mathcal{I}$.

\begin{definition}
A symmetric $(v,k)$-configuration
$\mathcal{I}=(\mathcal{P},\mathcal{L})$ is \emph{cyclic} if
there exists a permutation $\sigma$ of $\mathcal{P}$ mapping
blocks to blocks, and acting regularly on both $\mathcal{P}$
and $ \mathcal{L}$.
\end{definition}

We recall two well-known examples from Finite Geometry.

\begin{example}
\label{SingerP} Any Desarguesian projective plane $PG(2,q)$ is a cyclic
symmetric $(q^{2}+q+1,q+1)$-configuration.
\end{example}

\begin{example}
\label{SingerA} Fix a point $P$ and a line $\ell $ in $PG(2,q)$
with $ P\notin \ell $. Let $\mathcal{P}$ be the point set
consisting of the points of $PG(2,q)$ distinct from $P$ and not
lying on $\ell $. Let $\mathcal{L}$ be the line set consisting
of lines of $PG(2,q)$ distinct from $\ell $ and not passing
through $P$. Then $(\mathcal{P},\mathcal{L})$ is a cyclic
symmetric $(q^{2}-1,q)$-configuration. In \cite{FunkLabNap}
such a configuration is called an \emph{anti-flag}.
\end{example}

For a cyclic symmetric $(v,k)$-configuration
$\mathcal{I}=(\mathcal{P}, \mathcal{L})$ let $S$ be the cyclic
group generated by $\sigma$. Let $ \mathcal{P}={P_0,\ldots,
P_{v-1}}$ and $\mathcal{L}={\ell_0,\ldots, \ell_{v-1}}$.
Arrange indexes in such a way that
\begin{eqnarray}
\sigma &:&P_{i}\mapsto P_{i+1\pmod v},  \notag \\
\ell _{i}&=&\sigma ^{i}(\ell _{0}).  \label{form4_li}
\end{eqnarray}
Clearly,
\begin{equation}
P_{i}=\sigma ^{i}(P_{0}),\quad P_{i+\delta }=\sigma ^{\delta }(P_{i}).
\label{form4_Pi}
\end{equation}

For any divisor $d$ of $v$ the group $S$ has an unique cyclic
subgroup $ \widehat{S}_{d}$ of order $d$ generated by $\sigma
^{t}$ where
\begin{equation}
t=\frac{v}{d}\text{.}  \label{form4_td}
\end{equation}

Let $O_{0},O_{1},\ldots ,O_{t-1}$ be the orbits of
$\mathcal{P}$ under the action of the
subgroup$~\widehat{S}_{d}$. Clearly, $|O_{i}|=d.$ We arrange
indexes so that $P_{0}\in O_{0}$ and $O_{w}=\sigma
^{w}(O_{0}),$ whence $ \sigma ^{w}(O_{c})=O_{c+w\pmod t}$.
Therefore,
\begin{equation}
O_{i}=\{P_{i},\sigma ^{t}(P_{i}),\sigma ^{2t}(P_{i}),\ldots ,\sigma
^{(d-1)t}(P_{i})\},
\text{ }i=0,1,\ldots ,t-1.  \label{form4_Oi}
\end{equation}
By (\ref{form4_Pi}),(\ref{form4_Oi}), each orbit $O_{i}$
consists of $d$ points $P_{u}$ with $u$ equal to $i$ modulo
$t$.

Let $L_{0},\ldots ,L_{t-1}$ be the orbits of the set $\mathcal{L}$ under the
action of $\widehat{S}_{d}$. Obviously, $|L_{i}|=d.$ We arrange indexes in
such a way that $\ell _{0}\in L_{0}$ and $L_{w}=\sigma ^{w}(L_{0})$. As a
result,
\begin{equation}
L_{i}=\{\ell _{i},\sigma ^{t}(\ell _{i}),\sigma ^{2t}(\ell _{i}),\ldots
,\sigma ^{(d-1)t}(\ell _{i})\},\text{ }i=0,1,\ldots ,t-1.  \label{form4_Li}
\end{equation}

Let
\begin{equation}
w_{u}=|\ell _{0}\cap O_{u}|,\text{ }u=0,1,\ldots ,t-1.  \label{form4_wu}
\end{equation}
Clearly,
\begin{equation}
w_{0}+w_{2}+\ldots +w_{t-1}=k.  \label{form4_sum_wi}
\end{equation}

\begin{theorem}
\label{Th4_the-same-meeting} Let
$\mathcal{I}=(\mathcal{P},\mathcal{L})$ be a cyclic symmetric
$(v,k)$-configuration. Let $d$ and $t$ be integers as in
\textrm{(\ref{form4_td})}. Let $O_{0},\ldots ,O_{t-1}$ and
$L_{0},\ldots ,L_{t-1}$ be$,$ respectively, orbits of points
and lines of $\mathcal{I}$ under the action of
$\widehat{S}_{d}$. Assume also that for points, lines, and
orbits, indexes are arranged as in \textrm{(\ref{form4_li}),
(\ref {form4_Pi})}, \textrm{(\ref{form4_Oi}),
(\ref{form4_Li})}. Then for any $i$ and $j,$ every line of the
orbit $L_{i}$ meets the orbit $O_{j}$ in the same number of
points $w_{j-i\pmod t}$, where $w_{u}$ is defined by
\textrm{(\ref {form4_wu})}.
\end{theorem}

\begin{pro_f}
Fix some $j\in \{0,1,\ldots ,t-1\}$. For any $a=0,1,\ldots ,v-1$, let $s_{a}$
be such that $0\leq s_{a}\leq t-1$ and $s_{a}\equiv j-a\pmod t$. As $\sigma
^{a}$ maps the orbit $O_{s_{a}}$ on $O_{j}$,
\begin{equation*}
\mid \sigma ^{a}(l_{0})\cap O_{j}\mid =|l_{0}\cap O_{s_{a}}|=w_{s_{a}}
\end{equation*}
holds. This proves that for any $i=0,1,\ldots ,t-1$ the line
orbit $L_{i}$ consists of lines meeting $O_{j}$ in the same
number of points $w_{j-i\pmod t }$. Then, as $O_{j}$ and
$L_{i}$ have the same size, through any point $P\in O_{j}$
there pass exactly $w_{j-i\pmod t}$ lines in $L_{i}$, each of
which meets $O_{j}$ in $w_{j-i\pmod t}$ points.$\hfill \qed$
\end{pro_f}

\begin{definition}
An $i\times j$ matrix $A$ is said to be \emph{$d$-block circulant} if
\begin{equation}  \label{form4_incid-matr-planebis}
A=\left[ \renewcommand{\arraystretch}{0.85} \setlength{\arraycolsep}{5pt}
\begin{array}{cccc}
C_{0,0} & C_{0,1} & \ldots & C_{0,j/d-1} \\
C_{1,0} & C_{1,1} & \ldots & C_{1,j/d-1} \\
\vdots & \vdots & \vdots & \vdots \\
C_{i/d-1,0} & C_{i/d-1,1} & \ldots & C_{i/d-1,j/d-1}
\end{array}
\right]
\end{equation}
for some binary \emph{circulant }$d\times d$ matrices
$C_{\lambda,\mu}$, $ 0\le \lambda\le i/d-1,\, 0\le \mu\le
j/d-1$.

For a $d$-block circulant $i\times j$ matrix $A$ as in (\ref
{form4_incid-matr-planebis}), the weight matrix $W(A)$ of $A$
is the $\frac{i }{d}\times \frac{j}{d}$ matrix whose entry in
position $\lambda ,\mu $ is the weight of $C_{\lambda ,\mu }$,
that is, the number of units in each row of $C_{\lambda ,\mu
}$.
\end{definition}

We remark that in \cite{FunkLabNap} square $i\times i$ $d$-block circulant
matrices are called $(i,d)$-polycir\-culant.

\begin{corollary}
\label{Cor4_SingSubgPart} Let $\mathcal{I}$ be as in Theorem
\textrm{\ref {Th4_the-same-meeting}}. Then the $J_{4}$-free
incidence $v\times v$ matrix $ V$ of $\mathcal{I}$ is $d$-block
circulant. More precisely, it can be represented as a $t\times
t$ block matrix
\begin{equation}
V=\left[ \renewcommand{\arraystretch}{0.85} \setlength{\arraycolsep}{5pt}
\begin{array}{ccccc}
C_{0,0} & C_{0,1} & C_{0,2} & \ldots & C_{0,t-1} \\
C_{1,0} & C_{1,1} & C_{1,2} & \ldots & C_{1,t-1} \\
\vdots & \vdots & \vdots & \vdots & \vdots \\
C_{t-1,0} & C_{t-1,1} & C_{t-1,2} & \ldots & C_{t-1,t-1}
\end{array}
\right]  \label{form4_incid-matr-plane}
\end{equation}
where $C_{i,j}$ is a $J_{4}$-free binary \emph{circulant} $d\times d$ matrix
of weight $w_{j-i\pmod t}.$ The weight matrix $W(V)$ is the circulant matrix
\begin{equation}
W(V)=\left[ \renewcommand{\arraystretch}{0.80} \setlength{
\arraycolsep}{0.4em}
\begin{array}{ccccccc}
w_{0} & w_{1} & w_{2} & w_{3} & \ldots & w_{t-2} & w_{t-1} \\
w_{t-1} & w_{0} & w_{1} & w_{2} & \ldots & w_{t-3} & w_{t-2} \\
w_{t-2} & w_{t-1} & w_{0} & w_{1} & \ldots & w_{t-4} & w_{t-3} \\
w_{t-3} & w_{t-2} & w_{t-1} & w_{0} & \ldots & w_{t-5} & w_{t-4} \\
\vdots & \vdots & \vdots & \vdots & \vdots & \vdots & \vdots \\
w_{1} & w_{2} & w_{3} & w_{4} & \ldots & w_{t-1} & w_{0}
\end{array}
\right] .  \label{form4_circulweights}
\end{equation}
\end{corollary}

\begin{pro_f}
A submatrix $C_{i,j}$ corresponds to the intersections of the lines of the
orbit $L_{i}$ with the points of orbit $O_{j}.$ The submatrix is circulant
due to the arrangements of \textrm{(\ref{form4_li}), (\ref{form4_Pi})},
\textrm{(\ref{form4_Oi}), (\ref{form4_Li})}.$\hfill \qed$
\end{pro_f}

\begin{remark}
\label{Rem4_circ-weight-h} Let $C$ be a $J_{4}$-free binary
circulant $ d\times d$ matrix of weight $w$. Then several
$J_{4}$-free matrices with different parameters can be
constructed by using the natural \emph{\ decomposition of
square circulant }01-\emph{matrices}, cf. \cite[Sect\-ion
IV,~B]{DecompMask}. Such new matrices are circulant or consist
of circulant submatrices.

From now on, we assume that in circulant matrices rows are
shifted to the right. For a circulant matrix $C$, consider the
set $s(C)=\{s_{1},s_{2}, \ldots ,s_{w}\}$ of the positions of
the units in the first row of $C,$ arranged in such a way that
$1\leq s_{1}<s_{2}<\ldots <s_{w}\leq d$.

Let $I_{d}=I_{d}(1)$ be the identity matrix of order $d$ and let $I_{d}(v)$
be the circulant permutation $d\times d$ matrix obtained from $I_{d}$ by
shifting of every row by $v-1$ positions. The matrix $C$ can be viewed as
the superposition of $w$ matrices $I_{d}(s_{i}),$ $i=1,2,\ldots ,w.$

Taking into account that the \textquotedblleft starting\textquotedblright\
matrix $C$ is $J_{4}$-free, it is easy to see that the following holds:

\begin{enumerate}
\item The superposition of any $w-\delta $ distinct matrices $I_{d}(s_{i})$
gives a $J_{4}$-free circulant $d\times d$ matrix of weight $w-\delta$.

\item Let $S_{i}$ be a subset of $s(C)$ with
    $|S_{i}|=w-\delta_{i}$. Let $ C(S_{i})$ be the
    superposition of $w-\delta_{i}$ matrices $I_{d}(s_{u})$
with $s_{u}\in S_{i}$. A matrix consisting of distinct
submatrices of type $ C(S_{i})$ is $J_{4}$-free provided
that the subsets $S_{i}$ are pairwise disjoint.

\item If the subsets $S_{1},S_{2}$ are disjoint, then the following matrix
is $J_{4}$-free:
\begin{equation*}
\left[ \renewcommand{\arraystretch}{1.0}
\begin{array}{cc}
C(S_{1}) & C(S_{2}) \\
C(S_{2}) & C(S_{1})
\end{array}
\right] .
\end{equation*}
\end{enumerate}
\end{remark}

\begin{remark}
\label{Rem4_circsubm} Let a matrix $M(m,n)$ be $d$-block circulant. Then, by
Remark \ref{Rem4_circ-weight-h}, a matrix $M(m,n-\delta )$ can be
constructed for any $\delta \leq n$.
\end{remark}

We are now in a position to describe our second construction of matrices
\newline
$M(m_{1},m_{2},n_{2},n_{2})$. From now on, the subscript difference $j-i$ is
calculated modulo $t$. \medskip

\textbf{Construction B.} Let
$\mathcal{I}=(\mathcal{P},\mathcal{L})$ be a cyclic symmetric
$(v,k)$-configuration. Let $d$ and $t$ be integers as in
\textrm{(\ref{form4_td})}. Let $V$ be as
in~(\ref{form4_incid-matr-plane}). Fix some non-negative
integers $u_{1},\ldots ,u_{r}$, $0\leq u_{i}\leq t-1$. Let
$V^{\prime }$ be a matrix obtained from $V$ by replacing the
circulant submatrices $C_{i,j}$ such that $j-i=u_{k}$ with
$d\times d$ matrices $ C_{i,j}(S_{k})$ with
$|S_{k}|=w_{j-i}-\delta _{k}$, as in Remark~\ref
{Rem4_circ-weight-h}. Let $W(V^{\prime })$ be the weight matrix
of $ V^{\prime }$. If an $\frac{m_{1}}{d}\times
\frac{m_{2}}{d}$ submatrix of $ W(V^{\prime })$ is such that
the sum of the elements of every row (column) is equal to
$n_{1}$ ($n_{2}$), then the corresponding submatrix of $
V^{\prime }$ is a $J_{4}$-free matrix
$M(m_{1},m_{2},n_{1},n_{2})$. For $ m_{1}=m_{2}$ and
$n_{1}=n_{2},$ a matrix $M(m,n)$ is obtained.

\begin{example}
\label{Ex4_circul} The matrices $C_{i,j}{(S)}$ obtained from
the submatrix $ C_{i,j}$ of (\ref{form4_incid-matr-plane}) as
in Remark~\ref {Rem4_circ-weight-h}, are circulant matrices
$M(d,w_{j-i}-\delta )$, where $ |S|=w_{j-i}-\delta $.
Therefore, a family of $J_{4}$-free circulant matrices with the
following parameters is obtained (cf.
Remark~\ref{Rem4_circsubm}):
\begin{equation*}
M(m,n){}:{}m=d,\text{ }n=w_{u}-\delta ,\text{ }u=0,1,\ldots ,t-1,\,\,\delta
=0,1\ldots ,w_{u}-1.
\end{equation*}
\end{example}

\begin{example}
\label{Ex4} By Remark \ref{Rem4_circ-weight-h}, for any
sequence of non-negative integers $w_{0}^{\prime
},w_{1}^{\prime },\ldots ,w_{t-1}^{\prime }$ such that
$w_{i}^{\prime }\leq w_{i}$ there exists a $d$ -block circulant
submatrix $V^{\prime }$ of $V$ with
\begin{equation*}
W(V^{\prime })=\left[ \renewcommand{\arraystretch}{0.80}\setlength{
\arraycolsep}{0.4em}
\begin{array}{ccccccc}
w_{0}^{\prime } & w_{1}^{\prime } & w_{2}^{\prime } & w_{3}^{\prime } &
\ldots & w_{t-2}^{\prime } & w_{t-1}^{\prime } \\
w_{t-1}^{\prime } & w_{0}^{\prime } & w_{1}^{\prime } & w_{2}^{\prime } &
\ldots & w_{t-3}^{\prime } & w_{t-2}^{\prime } \\
w_{t-2}^{\prime } & w_{t-1}^{\prime } & w_{0}^{\prime } & w_{1}^{\prime } &
\ldots & w_{t-4}^{\prime } & w_{t-3}^{\prime } \\
w_{t-3}^{\prime } & w_{t-2}^{\prime } & w_{t-1}^{\prime } & w_{0}^{\prime }
& \ldots & w_{t-5}^{\prime } & w_{t-4}^{\prime } \\
\vdots & \vdots & \vdots & \vdots & \vdots & \vdots & \vdots \\
w_{1}^{\prime } & w_{2}^{\prime } & w_{3}^{\prime } & w_{4}^{\prime } &
\ldots & w_{t-1}^{\prime } & w_{0}^{\prime }
\end{array}
\right] .  \label{form4_circulweightsbis}
\end{equation*}
By Construction B, several $J_{4}$-free matrices $M(m,n)$ and $
M(m_{1},m_{2},n_{1},n_{2})$ can be obtained as submatrices of
$V^{\prime }$. Here, we provide a list of parameters
$m,n,m_{1},m_{2},n_{1},n_{2}$ for some of these matrices.
Significantly, every such matrix consists of circulant
submatrices$.$

\begin{enumerate}
\item The matrix $V^{\prime }$ itself is a matrix in $M(m,n)$. Taking into
account (\ref{form4_sum_wi}), the possible choices for $V^{\prime }$ give
rise to matrices of type
\begin{equation*}
M(m,n):m=v,\,\,n=\sum_{u=0}^{t-1}(w_{u}-\delta _{u})=k-\delta ,\text{ }
\delta \leq k-1.
\end{equation*}

\item Assume that $w_{0}^{\prime }=w_{1}^{\prime }=\ldots
    =w_{u-1}^{\prime }=w$, $u\geq 2$. Then for each
    non-negative integer $c\leq \lceil \frac{u}{2} \rceil $
    the submatrix of $V^{\prime }$
\begin{equation*}
(C_{i,j}^{\prime })_{i=c-1,\ldots ,2c-2,\,\,j=0,\ldots ,c-1}
\end{equation*}
is a matrix of type
\begin{equation*}
M(m,n):m=cd,\text{ }n=cw,\,c=1,2,\ldots ,\left\lceil \frac{u}{2}\right\rceil
.
\end{equation*}
Matrices of type
\begin{equation*}
M(m,n):m=cd,\text{ }n=(c-h)w,\,\,c=1,2,\ldots ,\left\lceil \frac{u}{2}
\right\rceil ,\,h=1,\ldots ,c-1.
\end{equation*}
can be easily obtained by substituting some of the
submatrices $ C_{i,j}^{\prime }$ with the null $d\times d$
matrix. Finally, the submatrices of $V^{\prime }$
\begin{equation*}
(C_{i,j}^{\prime })_{i=c_{1}-1,\ldots ,c_{1}+c_{2}-2,\,\,j=0,\ldots
,c_{1}-1},\text{ }1\leq c_{1}\leq u,~1\leq c_{2}\leq u+1-c_{1}
\end{equation*}
are of type
\begin{eqnarray*}
M(m_{1},m_{2},n_{1},n_{2}) &:&m_{1}=c_{1}d,\text{ }m_{2}=c_{2}d,\text{ }
n_{1}=c_{2}w, \\
n_{2} &=&c_{1}w,\text{ }c_{1}=1,2,\ldots ,u,\text{ }c_{2}=1,2,\ldots
,u+1-c_{1}.
\end{eqnarray*}

\item Similarly to ii), if $w_{1}^{\prime }=\ldots =w_{t-1}^{\prime
}=w,w_{0}^{\prime }\neq w$ then matrices of the following types can be
obtained:
\begin{equation*}
M(m,n):m=cd,~n=w_{0}^{\prime }+(c-h)w,\,\,c=2,3,\ldots ,t-1,~h=1,2,\ldots ,c.
\end{equation*}

\item Assume that $w_{i}^{\prime }=w_{i+f}^{\prime }=w_{i+f+u}^{\prime
}=w_{i+2f+u}^{\prime }=w,f,u\geq 1$. Then the submatrix of $V^{\prime }$
\begin{equation*}
\left[
\begin{array}{cc}
C_{0,i+f}^{\prime } & C_{0,i+2f+u}^{\prime } \\
C_{f,i+f}^{\prime } & C_{f,i+2f+u}^{\prime }
\end{array}
\right]
\end{equation*}
is a matrix of type
\begin{equation*}
M(m,n):m=2d,\text{ }n=2w.
\end{equation*}

\item Assume that $w_{u+1}^{\prime }=w_{0}^{\prime },$ $w_{u+2}^{\prime
}=w_{1}^{\prime },\ldots ,w_{2u}^{\prime }=w_{u-1}^{\prime },$ $u\geq 1$.
Then the submatrix of $V^{\prime}$
\begin{equation*}
(C^{\prime}_{i,j})_{i=u,\ldots,2u,\,\,j=0,\ldots,u}
\end{equation*}
is a matrix of type
\begin{equation*}
M(m,n):m=(u+1)d,\text{ }n=w_{0}^{\prime }+w_{1}^{\prime }+...+w_{u}^{\prime
}.
\end{equation*}

\item Assume that $w_{i_{1}}^{\prime }=w_{i_{2}}^{\prime
    }=\ldots =w_{i_{u}}^{\prime }=w,$ $2\leq u\leq t.$ Then
    for any $c$-subset $B$ of $ \{i_{1},i_{2},\ldots
    ,i_{u}\}$ the matrix
\begin{equation*}
(C_{0,j}^{\prime })_{j\in B}
\end{equation*}
is a matrix of type
\begin{equation*}
M(m_{1},m_{2},n_{1},n_{2}) :m_{1}=d,\text{ }m_{2}=cd,\text{ }n_{1}=cw,
\text{ } n_{2}=w,\text{ }c=2,3,\ldots ,u.
\end{equation*}
\end{enumerate}
\end{example}

\begin{remark}
\label{Rem4_skeleton_circulant} Assume that $M(m,n)$ is
circulant (an instance is provided by
Example~\ref{Ex4_circul}). Let $ s(M(m,n))=(s_{1},s_{2},\ldots
,s_{n})$, cf. Remark \ref{Rem4_circ-weight-h}. We consider
$M(m,n)$ as a superposition of $n$ circulant permutation $
m\times m$ matrices $I_{m}(s_{i}),$ $i=1,\ldots ,n.$ Then the
skeleton $ S(m,n)$ of the parity-check matrix $H$ of any BG
code with supporting graph $ G(M(m,n))$ is the following matrix
consisting of circulant permutation $ m\times m$ submatrices:
\begin{equation}
S(m,n)=\left[ \renewcommand{\arraystretch}{1.0}
\begin{array}{cccc}
I_{m} & I_{m} & \ldots & I_{m} \\
I_{m}(s_{1}) & I_{m}(s_{2}) & \ldots & I_{m}(s_{n})
\end{array}
\right] .  \label{form4_S(m,n)}
\end{equation}
Such a structure of the skeleton matrix can be useful for code
implementation. Assume that for constituent $[n,k_{t}]$ codes
$\mathcal{C} _{t}$ we have $\mathcal{C}_{1}=\ldots
=\mathcal{C}_{m},$ $\mathcal{C} _{m+1}=\ldots
=\mathcal{C}_{2m},$ (cf. Section~\ref{Sec2}). Let $
r_{t}=n-k_{t}$. Let also $[c_{j,1}^{(t)}c_{j,2}^{(t)}\ldots
c_{j,r_{t}}^{(t)}]$ be the $j$th column of a parity-check
matrix $H_{t}$ of the code $\mathcal{C}_{t}$. We choose parity
check matrices for constituent codes in such a way that
$H_{1}=\ldots =H_{m},$ $H_{m+1}=\ldots =H_{2m}.$ Then the
parity-check matrix $H$ corresponding to the skeleton (\ref
{form4_S(m,n)}) has the form
\begin{equation*}
H=\left[ \renewcommand{\arraystretch}{1.0}
\begin{array}{cccc}
c_{1,1}^{(1)}I_{m} & c_{2,1}^{(1)}I_{m} & \ldots & c_{n,1}^{(1)}I_{m} \\
\vdots & \vdots & \vdots & \vdots \\
c_{1,r_{1}}^{(1)}I_{m}\medskip & c_{2,r_{1}}^{(1)}I_{m} & \ldots &
c_{n,r_{1}}^{(1)}I_{m} \\
c_{1,1}^{(m+1)}I_{m}(s_{1}) & c_{2,1}^{(m+1)}I_{m}(s_{2}) & \ldots &
c_{n,1}^{(m+1)}I_{m}(s_{n}) \\
\vdots & \vdots & \vdots & \vdots \\
c_{1,r_{m+1}}^{(m+1)}I_{m}(s_{1}) & c_{2,r_{m+1}}^{(m+1)}I_{m}(s_{2}) &
\ldots & c_{n,r_{m+1}}^{(m+1)}I_{m}(s_{n})
\end{array}
\right] .
\end{equation*}
The matrix $H$ consists of circulant submatrices and defines a
QC code, cf.
\cite{KouLinFos-RediscNew},\cite{Fos-circul},\cite{GabidISIT},\cite
{SurveyLivaSLLR}-\cite{QCFiniteField}. QC codes can be encoded
with the help of shift-registers with relatively small
complexity, see e.g. \cite{LinCost}, \cite{QCEncoder}.
\end{remark}

\section{\label{Sec5}Intersection Numbers of Orbits of Singer Subgroups}

In this section, the general results of Section \ref{Sec4} are
applied to the special case $\mathcal{I}=PG(2,q)$ when
$v=q^2+q+1$, $k=q+1$. We treat points of $PG(2,q)$ as nonzero
elements of$~F_{q^{3}}.$ Elements $a,b$ of $ F_{q^{3}}$
correspond to the same point if and only if $a=xb,$ $x\in
F_{q}.$ All points can be represented by the set $\{\alpha
^{0},\alpha ^{1},\alpha ^{2},\ldots ,\alpha ^{q^{2}+q}\}$ where
$\alpha $ is a a primitive element of $~F_{q^{3}}.$ Every class
of elements of $F_{q^{3}}$ corresponding to the same point has
one representative into the set. The point represented by $
\alpha ^{i}$ is denoted by $P_{i}$, i.e.,
$PG(2,q)=\{P_{0},P_{1},P_{2}, \ldots ,P_{q^{2}+q}\}.$ It is
well-known that the map
\begin{equation*}
\sigma :P_{i}\mapsto P_{i+1\pmod{q^{2}+q+1}}
\end{equation*}
is a projectivity of $PG(2,q)$ acting regularly on the set of points and on
the set of lines in $PG(2,q)$. The group $S$ generated by $\sigma$ is called
the Singer group of $PG(2,q)$, whereas groups ${\widehat{S}}_d$ as defined
in Section \ref{Sec4} are said to be Singer subgroups of $PG(2,q)$.

We investigate the possible cardinalities of the intersections
of a fixed line of $PG(2,q)$ with point orbits of Singer
subgroups, i.e., the values of $w_{u},$ see (\ref{form4_wu})
and Theorem \ref{Th4_the-same-meeting}. It is a relevant
problem for the purposes of this paper, as the parameters of
the new $J_{4}$-free matrices, obtained by Construction B,
depend on the features of the sequence
$w_{0},w_1,\ldots,w_{t-1}$, see (\ref{form4_sum_wi}) and
Examples \ref{Ex4_circul},\ref{Ex4}.

Let $q=p^{h}$ for some prime $p$. We use the notations of the previous
section. The map
\begin{equation*}
\tau :P_{i}\mapsto P_{ip\pmod{q^{2}+q+1}}
\end{equation*}
is a permutation on the points of $PG(2,q)$. Actually, it is
easy to see that $\tau $ is a collineation of $PG(2,q)$, that
is, $\tau$ maps collinear points onto collinear points. It is
well known (see e.g. \cite[Section 2.3.1] {dembowski1969}) that
under the action of a cyclic collineation group of a finite
projective plane, the point set and the line set have the same
cyclic structure. Therefore, as $\tau$ fixes $P_0$, at least
one line has to be left invariant by $\tau$. Arrange indexes in
such a way that this fixed line is $\ell_0$, cf.
(\ref{form4_li}), (\ref{form4_wu}). By (\ref{form4_Oi}), $ \tau
$ acts on the set of orbits $O_{0},\ldots ,O_{t-1}$ as follows:
$\tau (O_{i})=O_{pi\pmod t}.$ The orbit $O_{0}$ always is fixed
by $\tau .$

Denote by $s(p,t)$ the number of orbits of ${\mathbb{Z}}_{t}\setminus \{0\}$
under the action of the permutation group generated by the map $i\mapsto
p\cdot i.$ Then the following result clearly holds.

\begin{proposition}
The cyclic group generated by $\tau $ acts on the set $O_{1},\ldots ,O_{t-1}$
with $s(p,t)$ orbits.
\end{proposition}

\begin{proposition}
\label{Prop5_t-prime}Let $t$ be a prime. Then $s(p,t)$ divides
$t-1$ and $ s(p,t)$ is the least integer $i$ such that $p\equiv
\beta ^{i}\pmod t$ for some primitive element $\beta \in
{\mathbb{Z}}_{t}$.
\end{proposition}

\begin{pro_f}
Let $e$ be the order of $p\pmod t$ in the multiplicative group
of ${\mathbb{Z }}_{t}$. Then $s(p,t)=\frac{t-1}{e}$, and
$p\pmod t$ is the $s(p,t)$-th power of a primitive element in
${\mathbb{Z}}_{t}$.$\hfill \qed$
\end{pro_f}

\begin{remark}
As $q^{2}+q+1$ is odd, $t$ is odd too. Also, $t\neq 5,$ as it is
straightforward to check that $q^{2}+q+1$ is not divisible by $5$ for any
prime power$~q$.
\end{remark}

\begin{proposition}
\label{Prop5_s(p,t) bounds}The following holds$:$

\begin{enumerate}
\item If $t\neq 3,$ then $s(p,t)\leq \frac{t-1}{2}.$

\item If $t$ does not divide $p^{2}-1$, then $s(p,t)\leq \frac{t-1}{3}.$

\item $s(p,3)=2$ if $p\equiv 1\pmod 3,$ and $s(p,3)=1$ if $p\equiv 2\pmod 3.
$

\item $s(p,7)=2$ if $p\equiv 2$ or $4\pmod 7,$ and $s(p,7)=1$ if $p\equiv 3$
or $5\pmod 7$.
\end{enumerate}
\end{proposition}

\begin{pro_f}
Explicit values of $s(p,3)$ and $s(p,7)$ are straightforward.
Let $m(p,t)$ be the smallest size of an orbit of
${\mathbb{Z}}_{t}\setminus \{0\}$ under the action of the group
generated by the map $i\mapsto p\cdot i.$ If $ p\not\equiv
1\pmod t,$ then $m(p,t)\geq 2.$ Note that as
$(p^{2h}+p^{h}+1)/t$ must be an integer, the cases $p\equiv \pm
1\pmod t$ can occur only for $ t=3. $ Also, $m(p,t)=2$ holds if
and only if $p^{2}\equiv 1\pmod t.$ Otherwise, $m(p,t)\geq
3.\hfill \qed$
\end{pro_f}

Let $v(t)$ be the number of distinct values of the integers
$w_{u}$ of~(\ref {form4_wu}).

\begin{lemma}
\label{Lem5_wu=wpu_v(p,t)<=} For any non-negative integer $u\le
t-1$, values $w_{u}$ and $w_{pu\pmod t}$ coincide. Moreover,
$v(t)\leq s(p,t)+1$, and each value of $w_{u}$, with at most
one exception $w_{0}$ , occurs at least $ (t-1)/s(p,t)$ times.
\end{lemma}

\begin{pro_f}
By our previous assumption the line $\ell _{0}$ is fixed by
$\tau .$ Therefore $w_{u}=w_{pu\pmod t}$. The other assertions
are straightforward.$ \hfill \qed$
\end{pro_f}

The case $s(p,t)=1,$ $v(t)\leq 2,$ is investigated in~\cite{HAM-PEN}.

\begin{remark}
\label{Rem5}By Proposition \ref{Prop5_s(p,t) bounds} and Lemma
\ref {Lem5_wu=wpu_v(p,t)<=}, $v(3)\leq 3$ if $p\equiv 1\pmod
3.$ Actually, in \cite{DC2} it is proved that equality holds.
\end{remark}

\begin{lemma}
\label{Lem5_treuno} Let $w_{i}$ be as in \textrm{(\ref{form4_wu})}. Then
\begin{equation}
\sum_{i=0}^{t-1}w_{i}^{2}=\frac{q^{2}+(t+1)q+1}{t}.  \label{form5-sum-wi^2}
\end{equation}
\end{lemma}

\begin{pro_f}
We consider the point orbit $O_{0}$ and the line orbits $L_{i}$
under the Singer subgroup$~\widehat{S}_{d}.$ By Theorem
\ref{Th4_the-same-meeting} and Corollary
\ref{Cor4_SingSubgPart}, for any $i=0,\ldots ,t-1,$ through any
point of $O_{0}$ there pass exactly $w_{-i\pmod t}$ lines in
$L_{i}$, each of which meets $O_{0}$ in $w_{-i\pmod t}$ points.
Therefore, the points on $ O_{0}$ can be counted as follows:
\begin{equation*}
d=\frac{q^{2}+q+1}{t}=1+\sum_{i=0}^{t-1}w_{i}(w_{i}-1).
\end{equation*}
Together with (\ref{form4_sum_wi}), this proves the
assertion.$\hfill \qed$
\end{pro_f}

\begin{corollary}
\label{Cor5} Let $O_{i_{1}},\ldots ,O_{i_{s(p,t)}}$ be orbit representatives
of the action of $\tau $ on the orbits $O_{1},\ldots ,O_{t-1}$. If $t$ is
prime then
\begin{eqnarray*}
w_{0}+\frac{t-1}{s(p,t)}\sum_{j=1}^{s(p,t)}w_{i_{j}} &=&q+1; \\
w_{0}^{2}+\frac{t-1}{s(p,t)}\sum_{j=1}^{s(p,t)}w_{i_{j}}^{2} &=&\frac{
q^{2}+(t+1)q+1}{t}.
\end{eqnarray*}
\end{corollary}

\begin{pro_f}
The assertions follow from (\ref{form4_sum_wi}) and Lemmas~\ref
{Lem5_wu=wpu_v(p,t)<=},\ref{Lem5_treuno}.$\hfill \qed$
\end{pro_f}

\begin{proposition}
\label{Prop5_t=3}Let $t=3$, and let $i_{0},i_{1},i_{2}$ be such
that $ w_{i_{0}}\leq w_{i_{1}}\leq w_{i_{2}}$,
$\{i_{0},i_{1},i_{2}\}=\{0,1,2\}$. Let
$A=(q-\sqrt{q}+1)/3,\,B=(q+\sqrt{q}+1)/3,$
$C=(q-2\sqrt{q}+1)/3,\,D=(q+2 \sqrt{q}+1)/3$. Then the
following holds:

\begin{enumerate}
\item If $q$ is a non square or $q=p^{2m},$ $p\equiv 1\pmod
    3$, then $ w_{i_0}<A,$ $A<w_{i_1}<B,$ $w_{i_2}>B.$

\item If $q=p^{4m+2}$ and $p\equiv 2\pmod 3$, then $\sqrt{q}\equiv 2\pmod 3$
and $w_{i_0}=w_{i_1}=A,$ $w_{i_2}=D.$

\item If $q=p^{4m}$ and $p\equiv 2\pmod 3$, then $\sqrt{q}\equiv 1\pmod 3$
and $w_{i_0}=C,$ $w_{i_1}=w_{i_2}=B.$
\end{enumerate}
\end{proposition}

\begin{pro_f}
For the sake of simplicity assume that $i_0=i$. By
(\ref{form4_sum_wi}), ( \ref{form5-sum-wi^2}),
$w_{0}+w_{1}+w_{2}=q+1,
w_{0}w_{1}+w_{1}w_{2}+w_{2}w_{3}=\frac{1}{2}((q+1)^{2}-\frac{1}{3}
(q^{2}+4q+1))=\frac{1}{3}(q^{2}+q+1)$. Therefore,
$w_{0},w_{1},w_{2}$ are the roots of the cubic polynomial
\begin{equation*}
T^{3}-(q+1)T^{2}+\frac{q^{2}+q+1}{3}T-w_{0}w_{1}w_{2}.
\end{equation*}
The real function $t\mapsto t^{3}-(q+1)t^{2}+\frac{1}{3}
(q^{2}+q+1)t-w_{0}w_{1}w_{2}$ is increasing in $]-\infty ,A]$, decreasing in
$[A,B]$, increasing in $[B,+\infty \lbrack $.

Clearly, when $q$ is not a square, $w_{i}=\frac{1}{3}(q\pm
\sqrt{q}+1)$ cannot hold. Also, taking into account Remark
\ref{Rem5}, if $p\equiv 1 \pmod 3$ we have three distinct
values of $w_{u}.$ Finally, if $p\equiv 2\pmod 3,$ then, by
Proposition \ref{Prop5_s(p,t) bounds} and Lemma \ref
{Lem5_wu=wpu_v(p,t)<=}, $w_{0}<w_{1}<w_{2}$ cannot hold.$\hfill
\qed$
\end{pro_f}

\begin{proposition}
Let $t=3$, and let $i_0,i_1,i_2$ be as in Proposition~\ref{Prop5_t=3}. Then
for any $i,j\in \{0,1,2\}$, the difference $4q-(3(w_{i}+w_{j})-2(q+1))^{2}$
equals $3s^{2}$ for some integer $s$.
\end{proposition}

\begin{pro_f}
For the sake of simplicity assume that $i_{0}=i$. Assume also
that $i=0$, $ j=1$. Consider the integers $X=w_{0}+w_{1}$,
$Y=w_{0}w_{1}$. By (\ref {form5-sum-wi^2}), $X=(q+1)-w_{2},$
$Y=\frac{1}{3}(q^{2}+q+1)-w_{2}X,$ whence
$Y=X^{2}-(q+1)X+\frac{1}{3}(q^{2}+q+1).$ As the quadratic
polynomial in $T$
\begin{equation*}
T^{2}-XT+Y=T^{2}-XT+X^{2}-(q+1)X+\frac{q^{2}+q+1}{3}
\end{equation*}
has two integer roots $w_{0},w_{1}$, we obtain that
\begin{equation*}
X^{2}-4\left( X^{2}-(q+1)X+\frac{q^{2}+q+1}{3}\right) =
\frac{4q-(3X-2(q+1))^{2}}{3}
\end{equation*}
is the square of an integer. Then the assertion follows.$\hfill \qed$
\end{pro_f}

From the above Propositions, we obtain an easy algorithm to compute the
possible $w_{0},w_{1},w_{2}$: for any integer $s<\frac{2}{\sqrt{3}}\sqrt{q}$
compute the difference $\Delta =4q-3s^{2}$. For the square values of$~\Delta
$, the integers $\frac{1}{3}(\sqrt{\Delta }+2(q+1))$ are the only possible
values of sums $w_{i}+w_{j}$.

We now consider the case $t>3$.

\begin{proposition}
\label{30may} Let $q=p^{r}$ be a square. Let $t$ be a prime
divisor of $ q^{2}+q+1$. Assume that $p \pmod t$ is a generator
of the multiplicative group of ${\mathbb{Z}}_{t}$. (If $t$ is
not prime, than we have to assume that the permutation group
generated by the map $i\mapsto p\cdot i$ acts transitively on
${\mathbb{Z}}_{t}\setminus \{0\}$). Then $t$ divides either $
q+\sqrt{q}+1$ or $q-\sqrt{q}+1$, and
\begin{eqnarray*}
w_{1}=w_{2}=\ldots =w_{t-1}=\frac{q+1\pm \sqrt{q}}{t}, \\
w_{0}=\frac{ q+1\pm (1-t)\sqrt{q}}{t}.
\end{eqnarray*}
\end{proposition}

\begin{pro_f}
As $\tau ^{i}(O_{1})=O_{p^{i}\pmod t}$, the hypothesis of $p$
being a generator of the multiplicative group of
${\mathbb{Z}}_{t}$ ensures that $ \tau $ acts transitively on
the orbits $O_{1},\ldots ,O_{t-1}$. By our assumptions $\ell
_{0}$ is a line fixed by~$\tau $. Then clearly $
w_{1}=w_{2}=\ldots =w_{t-1}.$ Now (\ref{form4_sum_wi}) and
(\ref {form5-sum-wi^2}) give $w_{0}+(t-1)w_{1}=(q+1)$, $
2w_{0}w_{1}+(t-2)w_{1}^{2}=(q^{2}+q+1)/t$. Then $w_{1}$ is a
root of
\begin{equation*}
tT^{2}-2(q+1)T+\frac{q^{2}+q+1}{t}.
\end{equation*}
This proves the assertion.$\hfill \qed$
\end{pro_f}

We remark that the Proposition \ref{30may} is not empty: for
example, for $ q=81$, $t=7$ its hypothesis are satisfied. Some
other values for which this holds are as follows: $q=2^{8},$
$t=13;$ $q=5^{4},$ $t=7;$ $q=2^{12},$ $ t=19; $ $q=3^{8},$
$t=7;$ $q=2^{16},$ $t=13;$ $q=17^{4},$ $t=7.$

\begin{proposition}
\label{Prop5_BaerSub}Let $q$ be a square and $d=v(q+\sqrt{q}+1),$ $v\geq 1.$
Then $t=\frac{1}{v}(q-\sqrt{q}+1)$ and $w_{0}=\sqrt{q}+v,$ $w_{1}=\ldots
=w_{t-1}=v.$
\end{proposition}

\begin{pro_f}
Let $v=1.$ Then the point orbits $O_{0},O_{1},\ldots ,O_{t-1}$
are Baer subplanes of $PG(2,q)$. The partition of the points
$PG(2,q)$ given by the orbits of the map $\sigma ^{t}$ is a
partition into Baer subplanes \cite[ Section 4.3]{Hirs}. Every
line of $PG(2,q)$ meets precisely one subplane in $ \sqrt{q}+1$
points, and the remaining subplanes in one point. If $v>1,$
every orbit $O_{i}$ is an union of $v$ subplanes.$\hfill \qed$
\end{pro_f}

\begin{example}
\label{Ex5_BaerSub} By Proposition \ref{Prop5_BaerSub} for
$v=1$, Remarks \ref{Rem4_circ-weight-h}, \ref{Rem4_circsubm},
and Example \ref{Ex4}iii one can obtain a family of matrices
\begin{eqnarray}
M(m,n) &:&m=c(q+\sqrt{q}+1),\,\,n=\sqrt{q}+c-\delta , \\
&&c=2,3,\ldots ,q-\sqrt{q},\,\,\delta =0,1,\ldots ,\sqrt{q}+c-1.  \notag
\end{eqnarray}
\end{example}

\begin{proposition}
\label{Prop5_KestEb}Let $q$ be a square and $d=q-\sqrt{q}+1.$
Then $t=q+ \sqrt{q}+1$ and $w_u\le 2$ for every $u$. There are
$v_{j}$ values of $u$ for which $w_{u}=j$, where
$v_{0}=\frac{1}{2} (q+\sqrt{q}),$ $v_{1}=\sqrt{q} +1,$ and
$v_{2}=\frac{1}{2}(q-\sqrt{q}).$
\end{proposition}

\begin{pro_f}
Each of the point orbits $O_{0},O_{1},\ldots ,O_{t-1}$ are Kestenband-Ebert
complete arcs \cite[Section 4.3]{Hirs}. The numbers $v_{m}$ are provided by
straightforward computation.$\hfill \qed$
\end{pro_f}

\begin{proposition}
\label{Prop5_IncompCap}Let $d=3.$ Then
$t=\frac{1}{3}(q^{2}+q+1)$ and $ w_u\le 2$ for every $u$. There
are $v_{j}$ values of $u$ for which $w_{u}=j$ , where
$v_{0}=\frac{1}{3}(q^{2}-2q+1),$ $v_{1}=q-1,$ and $v_{2}=1.$
\end{proposition}

\begin{pro_f}
Each orbit $O_{u}$ consists of $3$ non-collinear points.$\hfill \qed$
\end{pro_f}

\begin{remark}
\label{Rem5_refefGH} Note that the construction of \cite[Theorem 3.8]{GH}
can be viewed as a particular case of the construction of the present
section, cf. Propositions \ref{Prop5_BaerSub} and \ref{Prop5_KestEb}.
\end{remark}

For some $q$ and $d$, sequences $w_0,\ldots,w_{t-1}$ are given
in Table~I. They are obtained from both the results of
Sections~\ref{Sec4},\ref {Sec5} and computer search. As a
matter of notation, an entry $s_{i}$ in a sequence indicates
that $s$ should be repeated $i$ times.

Table I only shows some of the partitions arising from Singer
subgroups. In fact, one can obtain many other results. For
example, we obtained by computer the following pairs $
(d,w_{0}):(991,15),(721,15),(1093,15),(817,21),(1261,21),(1519,21),$
$ (2107,21),\newline
(3997,36),(5419,36),(3169,39),(3571,39),(6487,39),(6487,52),(7651,52),$
$ (7651,57),\newline (19,4),(61,7),(93,8),(127,9),(217,12),$ $
(399,14),(469,16),(1261,25),(1387,27),\newline (1519,28),$
$(2107,31),(3997,43),(4921,48),(5419,49)$. All correspon- ding
matrices $\newline M(m,n)$ are circulant, see Example
\ref{Ex4_circul}. By Remark~\ref {Rem4_circ-weight-h}, several
matrices $M(d,w_{0}-\delta ),$ $\delta \geq 0$, can be obtained
as well.

\section{\label{Sec5.5}Orbits of Affine Singer Groups}

In this section, the general results of Section \ref{Sec4} are
applied to the cyclic symmetric $(q^2-1,q)$-configuration
described in Example~\ref {SingerA}. Fix the point $P=(1:0:0)$
and the line $\ell: X_0=0$ in $PG(2,q)$ . Let $\mathcal{P}$ be
the point set consisting of the points of $PG(2,q)$ distinct
from $P$ and not lying on $\ell$. Let $\mathcal{L}$ be the line
set consisting of lines of $PG(2,q)$ distinct from $\ell$ and
not passing through $P$. We treat points of $\mathcal{P}$ as
nonzero elements of $ F_{q^{2}}.$ All points can be represented
by the set $\{\alpha ^{0},\alpha ^{1},\alpha ^{2},\ldots
,\alpha ^{q^{2}-2}\}$ where $\alpha $ is a primitive element of
$F_{q^{2}}.$ The point represented by $\alpha ^{i}$ is denoted
by $P_{i}$, i.e., $\mathcal{P} =\{P_{0},P_{1},P_{2},\ldots
,P_{q^{2}-2}\}.$ It is well-known that the map
\begin{equation*}
\sigma :P_{i}\mapsto P_{i+1\pmod{q^{2}-1}}
\end{equation*}
\begin{center}
\textbf{TABLE I. Construction B. Theoretical and computer
results for the case $ \mathcal{I}=PG(2,q)$}
\begin{equation*}
\renewcommand{\arraystretch}{1.0}
\begin{array}{|@{\,}r|@{\,}r|@{\,}r|@{\,}l|@{}r|@{\,}r|@{\,}r|@{\,}l|}
\hline
q & d & t & w_{0},w_{1},\ldots ,w_{t-1} & q & d & t & w_{0},w_{1},\ldots
,w_{t-1} \\ \hline
4 & 7 & 3 & 1,1,3 & 53 & 409 & 7 & 12,9,9,5,9,5,5 \\ \hline
4 & 3 & 7 & 2,1,1,0,1,0,0 & 61 & 1261 & 3 & 16,25,21 \\ \hline
7 & 19 & 3 & 1,4,3 & 61 & 291 & 13 & 4_{2},3_{2},6_{2},4_{2},6,4,11,4,3 \\
\hline
7 & 3 & 19 & 2,0_{7},1,0,1,0,1,1,0,1,0,0,1 & 64 & 1387 & 3 & 27,19,19 \\
\hline
9 & 13 & 7 & 4,1,1,1,1,1,1 & 64 & 219 & 19 & 11,3_{18} \\ \hline
9 & 7 & 13 & 2,2,0,0,2,0,1,0,1,0,0,1,1 & 64 & 73 & 57 & 9,1_{56} \\ \hline
11 & 19 & 7 & 0,1,1,3,1,3,3 & 67 & 1519 & 3 & 21,28,19 \\ \hline
11 & 7 & 19 & \text{0,1,2,1}_{2}\text{,0,1,0}_{2}\text{,1,0,2,0}_{5}\text{
,2,1} & 67 & 651 & 7 & 8,13,13,7,13,7,7 \\ \hline
13 & 61 & 3 & 4,7,3 & 73 & 1801 & 3 & 19,28,27 \\ \hline
16 & 91 & 3 & 7,7,3 & 79 & 2107 & 3 & 31,21,28 \\ \hline
16 & 39 & 7 & 2,1,1,4,1,4,4 & 79 & 903 & 7 & 8,9,9,15,9,15,15 \\ \hline
16 & 21 & 13 & 5,1,1,1,1,1,1,1,1,1,1,1,1 & 81 & 949 & 7 & 4,13,13,13,13,13,13
\\ \hline
16 & 13 & 21 & \text{1,2}_{2}\text{,0}_{5}\text{,2,0,1,0,1,2,0}_{2}\text{,2,1
}_{2}\text{,0,2} & 81 & 511 & 13 & \text{4,8,8,6,9,8,9,1,6,4,6,9,4} \\ \hline
19 & 127 & 3 & 4,7,9 & 83 & 367 & 19 &
\begin{array}{l}
5,8,4,4,4,1,1,8,3,5, \\
4,4,5,1,5,5,8,5,4
\end{array}
\\ \hline
23 & 79 & 7 & 0,5,5,3,5,3,3 & 97 & 3169 & 3 & 28,31,39 \\ \hline
25 & 217 & 3 & 7,7,12 & 107 & 1651 & 7 & 13,24,15,15,13,15,13 \\ \hline
25 & 93 & 7 & 8,3,3,3,3,3,3 & 109 & 3997 & 3 & 31,43,36 \\ \hline
25 & 31 & 21 & 6,1_{20} & 109 & 1713 & 7 & 8,15,15,19,15,19,19 \\ \hline
25 & 21 & 31 &
\begin{array}{l}
2_{2},1_{2},2_{2},0,2,0_{2},1,0,2,1, \\
0,1,0_{3},1,2_{2},0_{3},2,0_{3},2,0
\end{array}
& 125 & 829 & 19 &
\begin{array}{l}
4,9,9,9,9,4,4,9,9,4, \\
9,9,4,4,4,4,9,4,9
\end{array}
\\ \hline
29 & 67 & 13 & 0,3,3,4,2,3,2,3,4,0,4,2,0 & 137 & 2701 & 7 &
15,23,15,23,23,24,15 \\ \hline
31 & 331 & 3 & 7,12,13 & 137 & 511 & 37 &
\begin{array}{l}
5,8_{2},2,5,4,5,2,6,2_{3},3, \\
5,2,4,3,2,5,2,4,5,3,2, \\
6,2_{2},3,8,3_{2},4_{3},2,0,6
\end{array}
\\ \hline
32 & 151 & 7 & 0,5,5,6,5,6,6 & 139 & 1497 & 13 & \text{10,16}_{2}\text{,10}
_{2}\text{,7,10}_{3}\text{,11,7,16,7} \\ \hline
37 & 469 & 3 & 13,9,16 & 149 & 3193 & 7 & 12,25,25,21,25,21,21 \\ \hline
37 & 201 & 7 & 8,3,3,7,3,7,7 & 151 & 3279 & 7 & 32,19,19,21,19,21,21 \\
\hline
37 & 67 & 21 &
\begin{array}{l}
0,1,2,0,1,2,1,1,2,4,1, \\
2,4,3,2,0,1,2,4,1,4
\end{array}
& 163 & 8911 & 3 & 63,49,52 \\ \hline
43 & 631 & 3 & 13,12,19 & 163 & 3819 & 7 & 32,25,25,19,25,19,19 \\ \hline
49 & 817 & 3 & 13,21,16 & 163 & 1407 & 19 &
\begin{array}{l}
6,12,6,6,9,6,7,6,9,7, \\
13,7,12,12,9,5,6,13,13
\end{array}
\\ \hline
49 & 57 & 43 & 8,1_{42} & 163 & 1273 & 21 &
\begin{array}{l}
9,9,5,7,9,7,7,7,5,11,5, \\
5,7,5,16,11,9,7,11,5,7
\end{array}
\\ \hline
\end{array}
\end{equation*}
\end{center}
is a collineation of the incidence structure
$\mathcal{I}=(\mathcal{P}, \mathcal{L})$ acting regularly on
both $\mathcal{P}$ and $\mathcal{L}$. The group $S$ generated
by $\sigma$ is called the affine Singer group of $
\mathcal{I}$, whereas groups ${\widehat{S}}_d$, as defined in
Section~\ref {Sec4}, are said to be affine Singer subgroups of
$\mathcal{I}$.

The cases $d=q+1$ and $d=q-1$ are of particular interest. Every
orbit under the action of the group ${\widehat{S}}_{q-1}$ is
precisely the intersection of a line of $PG(2,q)$ through $P$
with the point set $\mathcal{P}$. We will refer to this
intersection as to the \emph{trace} of the line. On the other
hand, the orbits under the action of the group
${\widehat{S}}_{q+1}$ are (disjoint) conics. Denote by $v(t)$
the number of distinct values of the integers $w_{u}$
of~(\ref{form4_wu}).

\begin{proposition}
\label{61} Let $t$ be a divisor of $q+1$. Then $v(t)=2$ and $
w_{1}=w_{2}=\ldots=w_{t-1}=(q+1)/t$, $w_{0}=(q+1)/t-1$.
\end{proposition}

\begin{pro_f}
Let $\ell ^{\prime }$ be a line of $\mathcal{L}$ and let $(\ell
^{\prime })^{\ast }$ be the line of $PG(2,q)$ through $P$
meeting the removed line $ \ell $ in the same point as~$\ell
^{\prime }$. Then $\ell ^{\prime }$ is disjoint from the trace
of $(\ell ^{\prime })^{\ast }$ but meets all the remaining $q$
traces. Each orbit $\mathcal{O}$ under the action of ${
\widehat{S}}_{d}$ with $d=(q^{2}-1)/t=(q-1)(q+1)/t$ is a union
of the traces of $(q+1)/t$ lines of $PG(2,q)$ through~$P$.
Therefore, the line $\ell ^{\prime }$ meets the orbit
$\mathcal{O}$ in $(q+1)/t$ points if $(\ell ^{\prime })^{\ast
}\notin \mathcal{O}$ and in $(q+1)/t-1$ points if $(\ell
^{\prime })^{\ast }\in \mathcal{O}$.$\hfill \qed$
\end{pro_f}

\begin{example}
\label{Ex6_t=q+1} Using Proposition \ref{61} for $t=q+1$,
Remarks \ref {Rem4_circ-weight-h}, \ref{Rem4_circsubm}, and
Examples \ref{Ex4}ii, \ref {Ex4}iii, one can obtain a family of
matrices
\begin{eqnarray}
M(m,n) &:&m=c(q-1),\,\,n=c-\delta ,\,\,\delta =0,1,\ldots ,c-1, \\
&&c=2,3,\ldots ,b,\,\,b=q\,\,\mathrm{if}\,\,\delta \geq 1,\,\,b=\left\lceil
\frac{q}{2}\right\rceil \,\,\mathrm{if}\,\,\delta =0.  \notag
\end{eqnarray}
\end{example}

\begin{proposition}
Let $t$ be a divisor of $q-1$. Then $v(t)\le 2(q-1)/t$ and
$w_{u}\le 2(q-1)/t $ for every $u$.
\end{proposition}

\begin{pro_f}
Each orbit under the action of ${\widehat{S}}_{d}$ with $
d=(q^{2}-1)/t=(q+1)(q-1)/t$ is a union of $(q-1)/t$ disjoint
conics. Therefore, each line of $\mathcal{L}$ meets this orbit
in at most $2(q-1)/t$ points.$\hfill \qed$
\end{pro_f}

For some $q$ and $t$ such that $t$ does not divide $q+1$,
sequences $ w_{0},\ldots ,w_{t-1}$ obtained from computer
search are given in Table II.

\section{\label{Sec6}Constructions not Arising from Collineation Groups}

\subsection{Product of Parabolas}

Let $AG(r,q)$ be the $r$-dimensional affine space over $F_{q}$
(which is sometimes denoted as the Euclidean space $EG(r,q).$)
A point in $AG(r,q)$ corresponds to a vector in $F_{q}^{r}$.
Following \cite{DGMP-even-caps}, in $ AG(2v,q)$ with $q$ even
we consider the product $K$ of $v$ parabolas, that is, the set
of size $q^{v}$
\begin{equation*}
K=\{(a_{1},a_{1}^{2},a_{2},a_{2}^{2},\ldots ,a_{v},a_{v}^{2})\mid
a_{1},\ldots ,a_{v}\in F_{q}\}\subset AG(2v,q).
\end{equation*}
In terms of \cite[Section 2]{Gi}, $K$ is a maximal translation
cap.

\begin{center}
\textbf{TABLE II. Construction B. Theoretical and computer
results for the case $ \mathcal{I}=AG(2,q)$}
\begin{equation*}
\renewcommand{\arraystretch}{1.0}
\begin{array}{|@{\,\,}r|@{\,\,}r|@{\,\,}r|@{\,\,}l|@{\,\,}r|@{\,\,}r|@{\,\,}r|@{\,\,}l|}
\hline
q & d & t & w_{0},w_{1},\ldots ,w_{t-1} & q & d & t & w_{0},w_{1},\ldots
,w_{t-1} \\ \hline
7 & 16 & 3 & 1,2,4 & 31 & 160 & 6 & 5,6,8,4,2,6 \\ \hline
7 & 8 & 6 & 1,2,2,0,0,2 & 37 & 456 & 3 & 9,12,16 \\ \hline
9 & 20 & 4 & 1,2,4,2 & 37 & 342 & 4 & 9,12,10,6 \\ \hline
9 & 10 & 8 & 1,2,2,2,0,0,2,0 & 37 & 114 & 12 & 1,4_{4},6,2,0,4,2,4,2 \\
\hline
11 & 24 & 5 & 3,2,2,0,4 & 41 & 420 & 4 & 13,8,8,12 \\ \hline
11 & 15 & 8 & 3,1,1,1,0,2,1,2 & 41 & 336 & 5 & 5,12,10,6,8 \\ \hline
11 & 12 & 10 & 1,2_{2},0,2_{2},0_{3},2 & 41 & 140 & 12 &
3,2_{3},5,2,4,5_{2},4,2,5 \\ \hline
13 & 56 & 3 & 5,2,6 & 43 & 616 & 3 & 17,10,16 \\ \hline
13 & 42 & 4 & 5,2,2,4 & 49 & 800 & 3 & 17,20,12 \\ \hline
13 & 28 & 6 & 3,0,2,2,2,4 & 49 & 600 & 4 & 9,12,16,12 \\ \hline
13 & 21 & 8 & 1,1,2,2,4,1,0,2 & 49 & 300 & 8 & 1,6,8,6,8,6,8,6 \\ \hline
13 & 14 & 12 & 1,0_{2},2_{4},0,2,0_{2},2 & 53 & 702 & 4 & 17,12,10,14 \\
\hline
17 & 72 & 4 & 5,2,4,6 & 61 & 1240 & 3 & 25,16,20 \\ \hline
17 & 36 & 8 & 1,0,2,2,4,2,2,4 & 61 & 930 & 4 & 13,18,18,12 \\ \hline
17 & 24 & 12 & 3,1,2_{2},1_{2},0,2,1,0,2_{2} & 81 & 1312 & 5 & 9,18,18,18,18
\\ \hline
19 & 120 & 3 & 9,6,4 & 81 & 656 & 10 & 1,8,10,8,10,8,10,8,10,8 \\ \hline
19 & 60 & 6 & 3,2,2,6,4,2 & 89 & 1980 & 4 & 25,18,20,26 \\ \hline
19 & 45 & 8 & 3,4,2,4,2,1,2,1 & 97 & 3136 & 3 & 37,34,26 \\ \hline
25 & 208 & 3 & 5,10,10 & 103 & 3536 & 3 & 41,32,30 \\ \hline
25 & 156 & 4 & 5,8,8,4 & 109 & 2970 & 4 & 29,32,26,22 \\ \hline
25 & 104 & 6 & 1,4,6,4,6,4 & 121 & 4880 & 3 & 33,44,44 \\ \hline
27 & 91 & 8 & 1,3,3,3,6,4,3,4 & 121 & 2440 & 6 & 13,20,24,20,24,20 \\ \hline
29 & 210 & 4 & 5,6,10,8 & 121 & 1220 & 12 &
1,10,12,10,12,10,12,10,12,10,12,10 \\ \hline
29 & 120 & 7 & 5,2,4,4,8,4,2 & 169 & 4080 & 7 & 13,26,26,26,26,26,26 \\
\hline
31 & 320 & 3 & 9,8,14 & 169 & 2040 & 14 &
\begin{array}{l}
1,12,14,12,14,12,14,12,14,12,14, \\
12,14,12
\end{array}
\\ \hline
31 & 192 & 5 & 7,4,6,4,10 &  &  &  &  \\ \hline
\end{array}
\end{equation*}
\end{center}

Let $A=AG(2v,q)\setminus K$ be the complement of $K$. Clearly,
$ |A|=q^{v}(q^{v}-1).$ In $AG(2v,q),$ there are
$q^{2v-1}B_{q,v}$ lines, where $B_{q,v}=(q^{2v}-1)/(q-1)$ is
the number of lines through every point. Since $K$ is a cap, it
is easily seen that $q^{v}B_{q,v}-q^{v}(q^{v}-1)/2$ lines meet
$K$ and $q^{v}(q^{v-1}-1)B_{q,v}+q^{v}(q^{v}-1)/2$ lines lie
entirely in $A.$

\begin{proposition}
\label{transcap} The number of lines contained in $A$ through a
given point is constant and equal to
$B_{q,v}-q^{v}+\frac{1}{2}(q-2).$
\end{proposition}

\begin{pro_f}
As $K$ is a maximal translation cap in $AG(2v,q),$ from
\cite[Proposition 2.5 ]{Gi} it follows that through any point
in $A$ there pass exactly $\frac{1}{2 }(q-2)$ secants of $K$.
In addition, we have exactly $q^{v}-q+2$ tangents to
$K$.$\hfill \qed$
\end{pro_f}

By Proposition \ref{transcap}, the incidence structure whose
points are the points in $A$ and whose lines are the lines
contained in $A$ is a configuration. Matrices in
$M(m_{1},m_{2},n_{1},n_{2})$ are then obtained for the
following values of the parameters:
\begin{eqnarray*}
m_{1} &=&q^{v}(q^{v-1}-1)\frac{q^{2v}-1}{q-1}+\frac{q^{v}(q^{v}-1)}{2},\text{
}n_{1}=q, \\
m_{2} &=&q^{v}(q^{v}-1),\text{ }n_{2}=\frac{q^{2v}-1}{q-1}-q^{v}+\frac{q-2}{2
},\text{ }q\text{ even.}
\end{eqnarray*}

\subsection{Projective Spaces and Subspaces}

Fix $h$ and $q$, and an integer $s$, with $0\leq s\leq h-1$. Consider the
following incidence structure: points are subspaces of $PG(h,q)$ of
dimension $s$; blocks are subspaces of $PG(h,q)$ of dimension $s+1$;
incidence is set-theoretical inclusion. This structure is a configuration.
By \cite[Theorem 3.1]{Hirs}, the numbers of points and blocks are,
respectively,
\begin{equation*}
v_{h,s}=\prod_{i=h-s+1}^{h+1}(q^{i}-1)\prod_{i=1}^{s+1}\frac{1}{q^{i}-1},
\text{ }b_{h,s}=\prod_{i=h-s}^{h+1}(q^{i}-1)\prod_{i=1}^{s+2}\frac{1}{q^{i}-1
};
\end{equation*}
the number of points in a block is $(q^{s+2}-1)/(q-1)$ (by
duality, the number of hyperplanes in a space of dimension
$s+1$ is the number of points of the space); the number of
blocks through a point is $(q^{h-s}-1)/(q-1)$ (again by
duality, the number of subspaces of dimension $s+1$ containing
a given subspace of dimension $s$ coincides with the number of
hyperplanes of a subspace of dimension $h-1-s$). Clearly, no
two points are contained in two distinct blocks. Then matrices
of the following type are obtained:
\setlength{\arraycolsep}{0.0em}
\begin{equation*}
M(m_{1},m_{2},n_{1},n_{2}){}:{}m_{1}=b_{h,s},\text{ }m_{2}=v_{h,s},\,\,n_{1}=
\frac{q^{s+2}-1}{q-1},\text{ }n_{2}=\frac{q^{h-s}-1}{q-1}.
\end{equation*}
When $h-s=s+2$, that is $h=2c$ is even and
$s=\frac{h}{2}-1=c-1$, we have $ b_{h,s}=v_{h,s},$
$n_{1}=n_{2}.$ This gives matrices
\begin{equation*}
M(m,n):m=\prod\limits_{i=c+2}^{2c+1}(q^{i}-1)\prod
\limits_{i=1}^{c}(q^{i}-1)^{-1},\text{ }n=\frac{q^{c+1}-1}{q-1},c\geq 1.
\end{equation*}
If $h=2$ the incidence structure is just the projective plane.
For $h=4$, $ c=2,$ $s=1,$ we obtain the following parameters
\begin{equation*}
M(m,n):m=(q^{2}+1)(q^{4}+q^{3}+q^{2}+q+1),n=q^{2}+q+1.
\end{equation*}

\subsection{$q$-Cancellation Construction}

This construction is given in \cite[Constructions 3.2,3.3]{GH}, see also the
references therein and \cite{AfDaZ}.

In the projective plane $PG(2,q)$ we fix a line $\ell $ and a
point $P$ and assign an integer $s\geq 0$. If $P\in \ell$ we
choose $s$ points on $\ell $ distinct from $P$, and $s$ lines
through $P$ distinct from $\ell$. If $ P\notin \ell$ we choose
$s$ arbitrary points on $\ell $ and consider the $s$ lines
connecting $P$ with these points. The incidence structure
obtained from $PG(2,q)$ by dismissing all the lines through the
$s+1$ selected points and all the points lying on the $s+1$
selected lines provides a matrix
\begin{equation}  \label{form7_q-canc}
M(m,n):n=q-s,\,\,m=\left\{
\begin{array}{@{\,\,}l@{\,\,}l}
q^{2}-qs & \mathrm{if}\,\,P\in \ell \\
q^2-(q-1)s-1 & \mathrm{if}\,\,P\notin \ell
\end{array}
\right..
\end{equation}

It should be noted that Example \ref{SingerA} describes a particular case of
this construction with $P\notin \ell $, $s=0$. Significantly, in this case
the incidence structure has a cyclic automorphism group. It essentially
extends the list of parameters and gives matrices consisting of circulant
submatrices, see Proposition~\ref{61} and Example~\ref{Ex6_t=q+1}.

\section{\label{Sec8}Summary of new symmetric configurations}

For $q$ a power prime, the constructions using complements of
Baer subplanes~ \cite{FunkLabNap} and Example~\ref{SingerP}
give rise to symmetric configurations with parameters
$(q^{4}-q,q^{2})$ and $(q^{2}+q+1,q+1)$. Such configurations,
together with configurations with parameters $(m,n)$ where $
m,n$ are as in~(\ref{form7_q-canc}), will be referred to as
\emph{classical} . In \cite{Gropp-nk} symmetric configurations
with the following parameters $ (m^{\ast },n^{\ast })$ are
obtained: $ (69,8),(89,9),(111,10),(145,11),(171,12),$\newline
$(213,13),$ $(255,14),(303,15),(355,16),(399,17),(433,18),$
$(493,19),$ $ (567,20),(667,21),\newline
(713,22),(745,23),(851,24),$ $(961,25)$. In \cite{Gropp-nk} it
is also proved that symmetric configurations with parameters
$(m_{1},n^{\ast })$ with $m_{1}>m^{\ast }$ exist for every pair
$(m^{\ast },n^{\ast })$. Therefore, it would be interesting to
obtain \emph{non-classical} configurations with parameters
$(m_{2},n^{\ast })$ with $m_{2}<m^{\ast }$.

In \cite{FunkLabNap} configurations with the following
parameters $ (m_{2},n^{\ast })$ are constructed:
$(98,10),\newline (242,14),(338,16),$ $(338,17),$
$(578,21),(722,23)$.

The following new parameters $(m_{2},n^{\ast })$ with
$m_{2}<m^{\ast }$ are obta\-ined from the constructions
proposed in Sections~\ref{Sec3}-\ref{Sec6} : $n^{\ast
}=8,\,m_{2}=65$; $n^{\ast }=11,\,m_{2}=133$; $n^{\ast
}=12,\,m_{2}=168$;\newline $n^{\ast }=13,\,m_{2}=183,189$;
$n^{\ast }=14,\,m_{2}=210,231,252$; $n^{\ast }=15,m_{2}=$
$231,252,272,273$; $n^{\ast
}=16,\,m_{2}=273,288,307,324,341,342 $; $n^{\ast
}=17,\,m_{2}=307,342,360,372,381$; $n^{\ast }=18$, \newline
$m_{2}=360,381,403$; $n^{\ast
}=19,\,m_{2}=381,434,462,465,484$; $n^{\ast }=20$, $m_{2}=465,$
$484,496,\newline 506,525,527,528,552,553,558,651$; $n^{\ast
}=21,\,m_{2}=496,506,527,552,553,558,576,589,\newline
598,600,620,624,644,650,651$; $n^{\ast
}=22,m_{2}=527,528,553,558,576,589,600,620,624,\newline
650,651,672,676,700,702$; $n^{\ast
}=23,\,m_{2}=553,558,589,600,620,650,651,672,676,700,\newline
702,728$; $n^{\ast }=24,$ $
m_{2}=589,620,624,651,676,702,728,756,780,784,806,810,812,837,\newline
840$; $n^{\ast
}=25,m_{2}=651,702,756,784,810,812,837,840,868,870,899,900,930
$.

Table III illustrates how to obtain some of the above
parameters $ m_{2} $. In column ``C'', $a,b$ and $c$ stand for
Examples \ref{Ex3_Hermit}, \ref{Ex5_BaerSub}, and
\ref{Ex6_t=q+1}, respectively, whereas $d$ (resp. $e$ ) means
that the parameters are obtained by applying Remark~\ref
{Rem4_circ-weight-h}i to the cyclic structures of Examples
\ref{SingerP} (resp. \ref{SingerA}); $f$ stands for Example
\ref{Ex4}iii.

\begin{center}
\textbf{TABLE III. Obtaining the new parameters}
\begin{equation*}
\renewcommand{\arraystretch}{1.0}
\begin{array}{|@{}c|@{\,}r|@{\,}c|@{\,}r|@{\,\,}r|@{\,\,}c||@{\,\,}c|@{\,\,}c|@{\,\,}c|@{\,\,}c|@{\,\,}c|@{\,\,}c|}
\hline
n^{\ast } & m_{2} & \text{C} & q & c & \delta & n^{\ast } & m_{2} & \text{C}
& q & c & \delta \\ \hline
8 & 65 & b & 9 & 5 & 0 & 16,17 & 342 & c & 19 & 19 & 3,2 \\ \hline
11 & 133 & d & 11 &  & 1 & 17,18 & 360 & d & 19 &  & 2,1 \\ \hline
12 & 168 & b & 16 & 8 & 0 & 17 & 372 & b & 25 & 12 & 0 \\ \hline
13 & 183 & d & 13 &  & 1 & \text{17,18,19} & 381 & d & 19 &  & \text{3,2,1}
\\ \hline
13 & 189 & b & 16 & 9 & 0 & 18 & 403 & b & 25 & 13 & 0 \\ \hline
14 & 210 & b & 16 & 10 & 0 & 19 & 434 & f & 25 & 2 &  \\ \hline
\text{14,15} & 231 & b & 16 & 11 & 1,0 & 19 & 462 & c & 23 & 21 & 2 \\ \hline
\text{14,15} & 252 & b & 16 & 12 & 2,1 & 19,20 & 465 & b & 25 & 15 & 1,0 \\
\hline
15 & 272 & c & 17 & 17 & 2 & 19,20 & 484 & c & 23 & 22 & 3,2 \\ \hline
\text{15,16} & 273 & d & 16 &  & 2,1 & 20,21 & 496 & b & 25 & 16 & 1,0 \\
\hline
16 & 288 & e & 17 &  & 1 & 20,21 & 506 & c & 23 & 23 & 3,2 \\ \hline
\text{16,17} & 307 & d & 17 &  & 2,1 & 20 & 525 & a & 25 &  &  \\ \hline
16 & 324 & c & 19 & 18 & 2 & \text{20,21,22} & 527 & b & 25 & 17 & \text{
2,1,0} \\ \hline
16 & 341 & b & 25 & 11 & 0 & 20 & 528 & c & 25 & 22 & 2 \\ \hline
\end{array}
\end{equation*}
\end{center}

\section*{Acknowledgments}

The authors would like to thank Valentine B. Afanassiev and
Victor V. Zyablov for useful discussions of code aspects of the
problems investigated, and Marien Abreu, Domenico Labbate and
Vito Napolitano for helpful discussions on symmetric
configurations.

\end{document}